\documentclass{birkjour_t2}
 \newtheorem{thm}{Theorem}[section]
 \newtheorem{cor}[thm]{Corollary}
 \newtheorem{lem}[thm]{Lemma}
 
 \theoremstyle{definition}
 \newtheorem{defn}[thm]{Definition}
 \theoremstyle{remark}
 \newtheorem{rem}[thm]{Remark}
 
 \numberwithin{equation}{section}
 \usepackage{cases}
 \usepackage[colorlinks=true,
 linkcolor=blue,
 anchorcolor=black,
 citecolor=blue, 
 urlcolor=blue
 ]{hyperref}
 \usepackage{cite}
 \usepackage{subfigure}
 \allowdisplaybreaks
 
\begin{document}

%-------------------------------------------------------------------------
% editorial commands: to be inserted by the editorial office
%
%\firstpage{1} \volume{228} \Copyrightyear{2004} \DOI{003-0001}
%
%
%\seriesextra{Just an add-on}
%\seriesextraline{This is the Concrete Title of this Book\br H.E. R and S.T.C. W, Eds.}
%
% for journals:
%
%\firstpage{1}
%\issuenumber{1}
%\Volumeandyear{1 (2004)}
%\Copyrightyear{2004}
%\DOI{003-xxxx-y}
%\Signet
%\commby{inhouse}
%\submitted{March 14, 2003}
%\received{March 16, 2000}
%\revised{June 1, 2000}
%\accepted{July 22, 2000}
%
%
%
%---------------------------------------------------------------------------
%Insert here the title, affiliations and abstract:
%

\title[A class of forward-backward regularizations...]
 {A class of forward-backward regularizations of the Perona-Malik equation with variable exponent}

%----------Author 1
\author[Y. Tong]{Yihui Tong}

\address{%
School of Mathematics\\
Harbin Institute of Technology\\
Harbin 150001 Heilongjiang\\
China}

\email{23B912019@stu.hit.edu.cn}

\author[W. Liu]{Wenjie Liu}

\address{%
School of Mathematics\\
Harbin Institute of Technology\\
Harbin 150001 Heilongjiang\\
China}

\email{liuwenjie@hit.edu.cn}

\author[Z. Guo]{Zhichang Guo}

\address{%
School of Mathematics\\
Harbin Institute of Technology\\
Harbin 150001 Heilongjiang\\
China}

\email{mathgzc@hit.edu.cn}

\author[W. Yao]{Wenjuan Yao}

\address{%
School of Mathematics\\
Harbin Institute of Technology\\
Harbin 150001 Heilongjiang\\
China}

\email{mathywj@hit.edu.cn}
%\thanks{This work was completed with the support of our
%\TeX-pert.}
%----------Author 2
%\author[J. Shao]{Jingfeng Shao}
%\address{The address of\br
%the second author\br
%sitting somewhere\br
%in the world}
%\email{mathywj@hit.edu.cn}

%----------classification, keywords, date


\date{January 1, 2004}
%----------additions
%\dedicatory{To my boss}
%%% ----------------------------------------------------------------------

\begin{abstract}
This paper investigates a novel class of regularizations of the Perona-Malik equation with variable exponents, of forward-backward parabolic type, which possess a variational structure and have potential applications in image processing. The existence of Young measure solutions to the Neumann initial-boundary value problem for the proposed equation is established via Sobolev approximation and the vanishing viscosity limit. The proofs rely on Rothe's method, variational principles, and Young measure theory. The theoretical results confirm numerical observations concerning the generic behavior of solutions with suitably chosen variable exponents.
\end{abstract}

%%% ----------------------------------------------------------------------
\maketitle
%%% ----------------------------------------------------------------------
%\tableofcontents
\section{Introduction}
In this paper, we consider the following forward-backward parabolic problem with nonstandard growth condition:
\begin{numcases}{}
	\displaystyle \frac{\partial u}{\partial t}=\operatorname{div}\left(\left(\frac{1}{1+\left|\nabla u\right|^2}+\delta\left|\nabla u\right|^{p(x)-2}\right)\nabla u\right), & $(x,t)\in Q_T$,\label{P1}\\
	\displaystyle \frac{\partial u}{\partial \vec{n}}=0, & $(x,t)\in \partial\Omega\times (0,T)$,\label{P2} \\
	u(x,0)=f(x), & $x\in\Omega$,\label{P3}
\end{numcases}
where $\Omega\subset\mathbb{R}^N$ ($N\geq 2$) is a bounded domain with smooth boundary $\partial\Omega$, $0<T<\infty$, $Q_T=\Omega\times(0,T)$, $\vec{n}$ denotes the unit outward normal to $\partial\Omega$, and $\delta>0$. The variable exponent function $p$ is globally $\log$-H\"{o}lder continuous and satisfies
\begin{equation}\label{p}
	2 \leq p^{-}=\operatorname*{ess\,inf}_{z\in \Omega} p(z)\leq p(x)\leq p^{+}=\operatorname*{ess\,sup}_{z\in \Omega} p(z)<\infty .
\end{equation}

\begin{figure}[htbp]
	\centering
	\includegraphics[scale = 0.70]{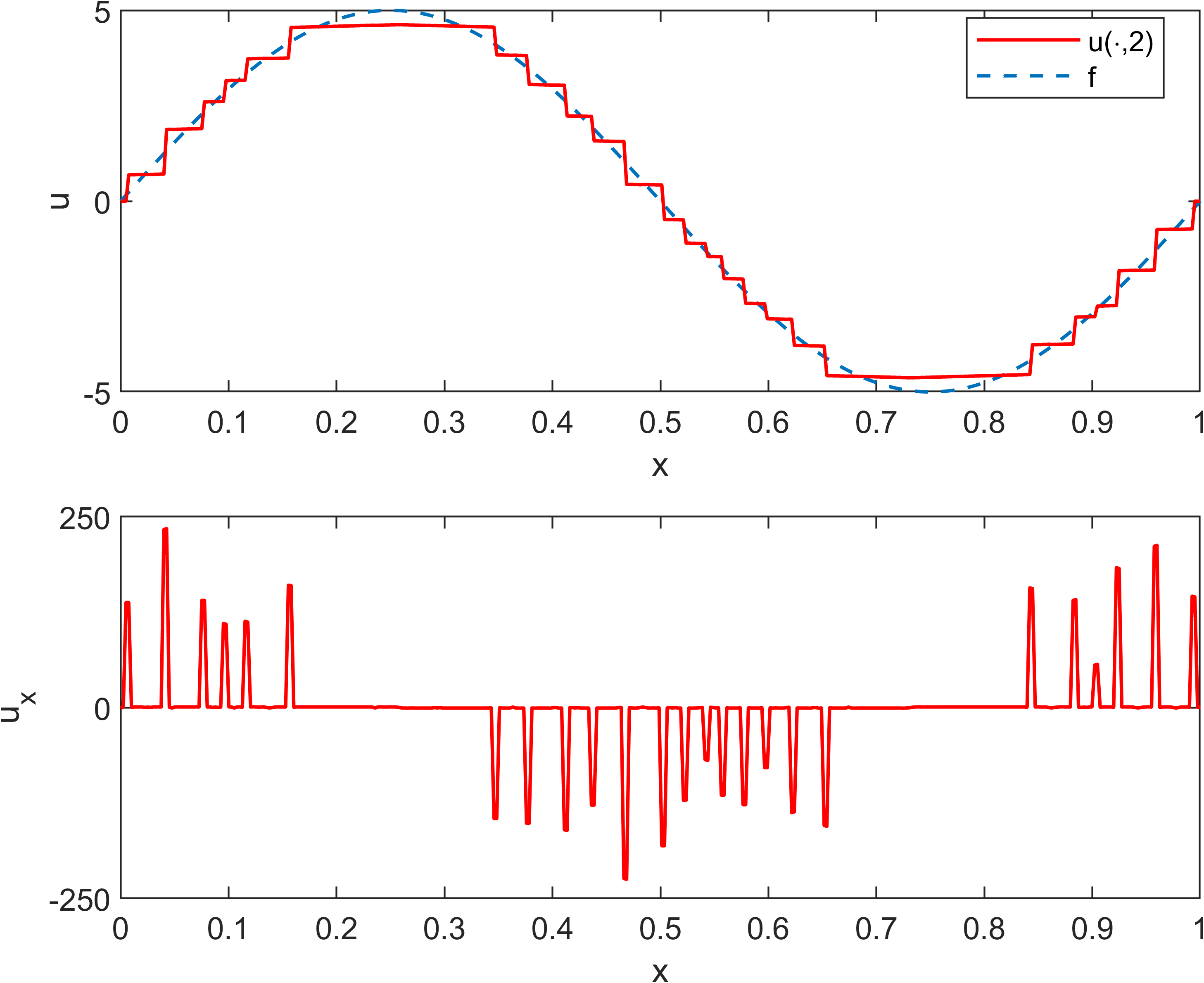}
	\caption{Manifestation of the staircasing phenomenon widely observed in discretizations of the Perona-Malik equation, where the graphs of $u(\cdot,t)$ and $u_x(\cdot,t)$ are plotted at time $t = 2$.}
	\label{PM fig}
\end{figure}

\begin{figure}[htbp]
	\centering
	\includegraphics[scale = 0.70]{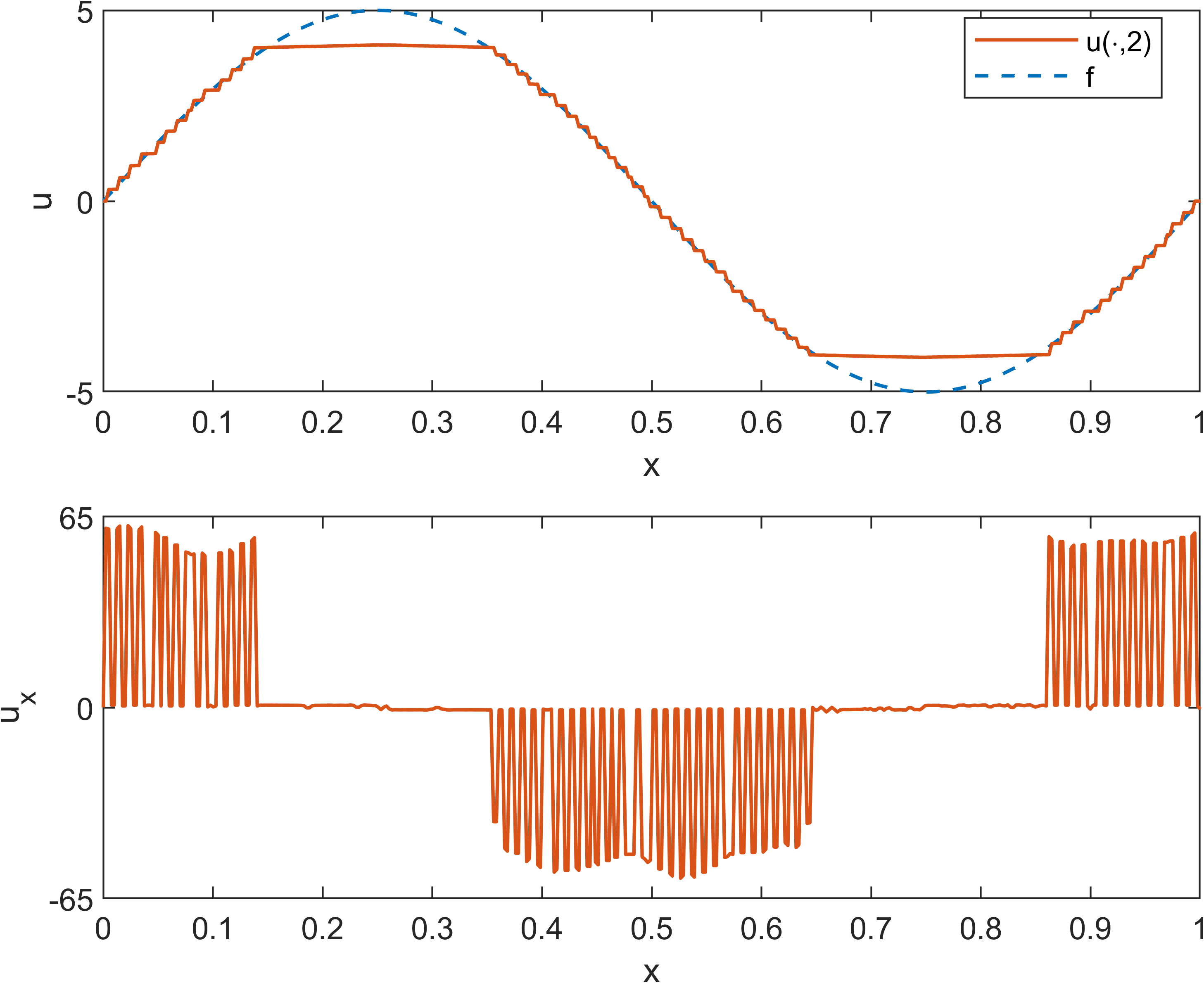}
	\caption{Microramping phenomenon observed for the regularized equation \eqref{P1} with constant exponent $p$ and $\delta = 0.001$.}
	\label{PLRPM fig}
\end{figure}

When $p>1$ is a constant, equation \eqref{P1} arises in \cite{Guidotti20131416} as a model for image denoising, where $u$ denotes the restored image, $f$ the noisy observation, and the function values quantify the gray levels. For sufficiently small $\delta$, this equation is of forward-backward parabolic type and admits no classical strong or distributional solutions. In \cite{Kinderlehrer1992001}, Young measure theory was employed to study the Dirichlet initial-boundary value problem (IBVP) for equation \eqref{P1} with $p=2$, where the authors proved the existence of Young measure solutions using Rothe's method and variational techniques combined with relaxation theory \cite{Dacorogna2008}. Demoulini \cite{Demoulini1996376} further investigated the qualitative properties of these solutions, including uniqueness, the comparison principle, and long-time behavior, which were later extended from the Dirichlet to the periodic and Neumann settings in \cite{GUIDOTTI20123226}. The well-posedness of the Dirichlet IBVP for equation \eqref{P1} was generalized to the case $1<p \leq 2$ by adopting an approach similar to that in \cite{Demoulini1996376}. For the case $p=1$, equation \eqref{P1} becomes singular, and Young measure solutions are considered in $L^{\infty}(0,T;BV(\Omega))$; see \cite{Yin2003659}. In fact, equation \eqref{P1} with constant exponent $p>1$ serves as a forward-backward regularization of the well-known Perona-Malik equation
\begin{equation}\label{PM}
	u_t=\operatorname{div} \left(\frac{\nabla u}{1+\left|\nabla u\right|^2}\right),
\end{equation}
which is parabolic if $\left|\nabla u\right|<1$ and backward-parabolic if $\left|\nabla u\right|>1$. Originally introduced for image reconstruction, the Perona-Malik equation is designed to reduce noise while preserving important image features, particularly edges and contours \cite{Perona1990629}. In practice, however, the Perona-Malik equation suffers from the drawback that it may result in staircase artifacts in the restored image. This issue inherently stems from the ill-posedness of \eqref{PM}, whose instabilities during the evolution manifest as the so-called staircasing phenomenon in numerical discretizations \cite{Guidotti2015131}; see Figure \ref{PM fig}. From a dynamical perspective, this phenomenon is replaced by a milder, microstructured ramping in the $p$-Laplacian regularized Perona-Malik equation, characterized by alternating gradients of finite amplitude; see Figure \ref{PLRPM fig}.

Forward-backward parabolic equations arise in various applications, including geophysics, biology, and the modeling of physical phenomena. A typical form can be written as
\begin{equation}\label{divFBform multidimen}
	u_t=\operatorname{div} \vec{\varPhi}(\nabla u),
\end{equation}
where the heat flux $\vec{\varPhi} (\xi) = \nabla \varPsi (\xi)$ for some nonconvex potential $\varPsi \in C^1(\mathbb{R}^N)$. In this case, the monotonicity condition $(\vec{\varPhi}(\xi)-\vec{\varPhi}(\eta)) \cdot (\xi-\eta) \geq 0$ is violated for some $\xi, \eta \in \mathbb{R}^N$. The one-dimensional version of equation \eqref{divFBform multidimen} reads
\begin{equation}\label{divFBform onedimen}
	u_t=\left[\varPhi(u_x)\right]_x,
\end{equation}
and appears independently as a model for heat and mass transfer in a stably stratified turbulent shear flow \cite{Barenblatt1993341}. Setting $w = u_x$, formal differentiation of \eqref{divFBform onedimen} with respect to $x$ leads to another type of forward-backward parabolic equation
\begin{equation}\label{DelFBform onedimen}
	w_t=\left[\varPhi(w)\right]_{xx},
\end{equation}
which arises in the theory of phase transitions \cite{Brokate2012}. In \cite{Hollig1983299}, the Neumann IBVP for \eqref{divFBform onedimen} is shown to admit infinitely many weak $L^2$-solutions when $\varPhi$ is a nonmonotonic piecewise linear function satisfying $s\varPhi(s)\geq C s^2$ for some $C>0$. This condition on $\varPhi$ was later relaxed in \cite{ZHANG2006322} by assuming higher regularity of the initial data. In \cite{Zhang2006171}, the existence of infinitely many weak $W^{1,\infty}$-solutions to the one-dimensional Perona-Malik equation under homogeneous Neumann boundary conditions was established. As shown in \cite{GURTIN1996178,Kichenassamy19971328}, no global $C^1$ solution exists for sufficiently large $\left|f^{\prime}\right|$, which suggests that even the local existence of solutions to the IBVP for equations \eqref{divFBform multidimen} or \eqref{divFBform onedimen} remains a nontrivial problem. Concerning the uniqueness of forward-backward parabolic problems, the situation becomes more subtle (see \cite{LAIR1985311}). For more theoretical results related to the Perona-Malik equation, the reader is referred to \cite{BOTTAZZI202139,CORLI2022474,ELLIOTT1985257,KIM20151889}.

A natural approach to handling ill-posed problems is to introduce suitable regularization, often motivated by physical considerations, which leads to well-posed problems \cite{Guidotti2015131}. Nevertheless, the qualitative behavior of solutions is highly sensitive to the type of regularization employed. In recent years, various regularization schemes have been investigated in different contexts (see \cite{BELLETTINI2008892,Giacomelli2019374,GUIDOTTI20094731,Shao20221779}). The pseudo-parabolic regularization of equation \eqref{divFBform onedimen},
\begin{equation}\label{pseudoRegu onedimen}
	u_t=\left[\varPhi(u_x)\right]_x + \varepsilon \left[\psi(u_x)\right]_{xt},
\end{equation}
with $\varepsilon > 0$ a small parameter and $\psi^{\prime} > 0$, has been proposed to account for time-delay effects \cite{Barenblatt1993341} and is applied to modeling the formation of layers with constant temperature or salinity in the ocean \cite{Barenblatt19931414}. By the transformation $w = u_x$, equation \eqref{pseudoRegu onedimen} is formally related to
\begin{equation}\label{SoboRegu onedimen}
	w_t=\left[\varPhi(w)\right]_{xx} + \varepsilon\left[\psi(w)\right]_{xxt},
\end{equation}
which arises in population dynamics \cite{Horstmann2004545}. Equation \eqref{SoboRegu onedimen} serves as the pseudo-parabolic regularization of equation \eqref{DelFBform onedimen}. If $\psi(u)=u$, this regularization is referred to as the Sobolev regularization \cite{Victor1998457}. The well-posedness of the Neumann IBVP for \eqref{pseudoRegu onedimen} was established in \cite{Barenblatt19931414}, where $\psi$ is an increasing but bounded function. It was also shown that solutions may develop spatial discontinuities (jumps) within finite time, with amplitudes that increase over time. This phenomenon results not only from the shape of $\varPhi$ but also from the boundedness of $\psi$, whereas it does not occur if $\psi$ has power-type or logarithmic growth at infinity \cite{Bertsch20131719,Bertsch2016051}. The existence of suitably defined nonnegative solutions to the Dirichlet IBVP for \eqref{SoboRegu onedimen} was proved in the space of Radon measures (see \cite{BERTSCH2015217,BERTSCH2016190,BERTSCH20172037,BERTSCH201846}), and results concerning uniqueness can be found in \cite{Smarrazzo2012085}. As $\varepsilon \to 0^{+}$, the family of solutions to the regularized problem is expected to converge to a measure-valued solution of the original IBVP for \eqref{DelFBform onedimen}; see \cite{Smarrazzo2013551}.

Both Dirichlet and Neumann IBVPs for the Cahn-Hilliard type regularization of equation \eqref{divFBform multidimen},
\begin{equation}\label{CHregu}
	u_t=\operatorname{div} \vec{\varPhi} (\nabla u) - \varepsilon \Delta^2 u,
\end{equation}
where $\vec{\varPhi}$ satisfies appropriate growth conditions, were investigated in \cite{Slemrod1989001} (see also \cite{Plotnikov1997550}). By the theory of Young measures (e.g., see \cite{ball2005version,Giaquinta1998,Valadier1990}), it has been proved that a subsequence of the family of approximating solutions to the regularized problem \eqref{CHregu} converges, in a suitable topology, to a Young measure solution of the ill-posed original IBVPs for \eqref{divFBform multidimen}. When the growth assumptions on $\vec{\varPhi}$ are relaxed, the Dirichlet IBVP for equation \eqref{divFBform multidimen} was studied in \cite{THANH20141403} through the limit $\varepsilon \rightarrow 0^{+}$ of solutions to the corresponding IBVP for Sobolev regularization
\begin{equation}\label{Sobolev regu}
	u_t=\operatorname{div} \vec{\varPhi} (\nabla u) + \varepsilon \Delta u_t ,
\end{equation}
where the term $\Delta u_t$ can be interpreted as accounting for viscous relaxation effects \cite{Novick1991331}. Specifically, the first step is to address equation \eqref{Sobolev regu} under the Dirichlet boundary condition for any $\varepsilon > 0$, and then to examine the vanishing viscosity limit $\varepsilon \rightarrow 0^{+}$ of the family $\{u^{\varepsilon}\}_{\varepsilon}$ of solutions to the approximating problems, obtaining Young measure solutions for the original problem. The asymptotic behavior of such solutions for large times is also studied in \cite{THANH20141403} using compactness arguments and $\omega$-limit set techniques. In view of the above results, a widely accepted strategy for handling forward-backward parabolic problems is to approximate them by a family of regularized problems, whose solutions are then used to extract limiting points, yielding suitably defined solutions of the original problem (see \cite{Evans20041599,Oliver20031380,Smarrazzo20101046,Smarrazzo2013551}). It should be noted that, except in specific cases (see \cite{Demoulini1996376,Mascia2009887}), the uniqueness of such solutions remains an open problem.

\begin{figure}[htbp]
	\centering
	\includegraphics[scale = 0.70]{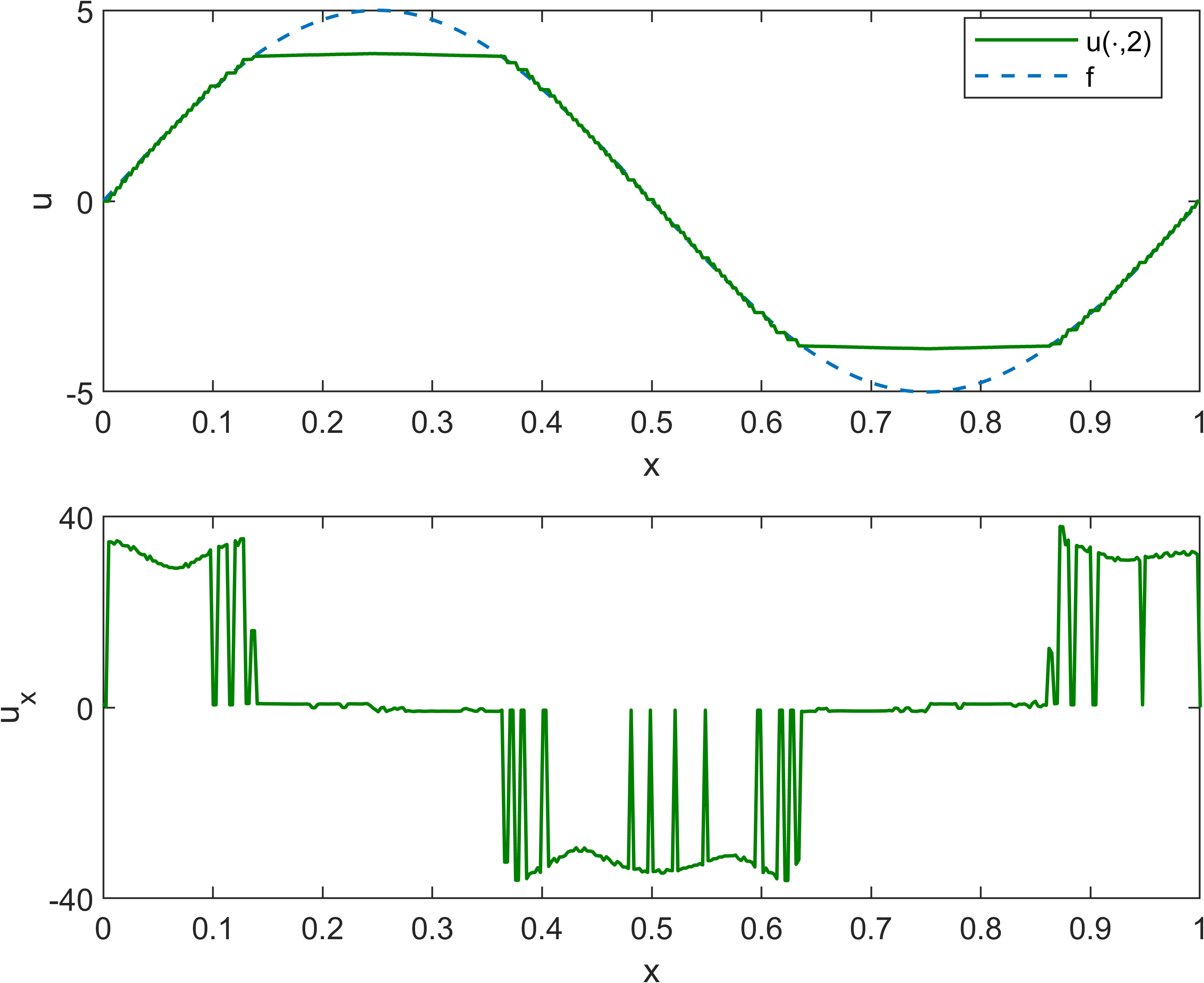}
	\caption{Gentler microramping phenomenon observed for the new regularized equation \eqref{P1} with $\delta = 0.001$.}
	\label{PxPlRPM fig}
\end{figure}

Here, we propose a novel class of regularizations \eqref{P1} for the Perona-Malik equation, interpretable as a $p(x)$-Laplacian regularization of forward-backward type, where the degree of regularization is adaptively modulated by spatial information. From the energy perspective, this regularization benefits from the variational structure
\begin{equation}\label{variational struc}
	\int_{\Omega} \left(\frac{1}{2}\log\left(1+\left|\nabla u\right|^2\right)+\frac{\delta}{p(x)} \left|\nabla u\right|^{p(x)}\right) dx ,
\end{equation}
which provides a mild variational regularization of the convex-concave energy functional underlying the Perona-Malik flow ($\delta=0$), while any $\delta>0$ is sufficient to make \eqref{variational struc} convex-concave-convex. Figures \ref{PM fig}--\ref{PxPlRPM fig} illustrate the transition from rough to smooth staircasing. In contrast to \eqref{P1} with a constant exponent, the solution to the proposed regularization for a suitable choice of $p$ shows better regularity; its gradients exhibit fewer oscillations, and the maximum gradients are confined to a smaller bounded interval. Meanwhile, equation \eqref{P1} can be regarded as an independent anisotropic diffusion model in image processing, where the spatially varying exponent $p$ is designed according to practical requirements. In this sense, the new model serves as an interpolation between two classical models: the Perona-Malik model and that introduced in \cite{CHEN20061383}. However, its application to image processing is not the primary focus of the present paper; interested readers are referred to \cite{CHEN20061383,Guo2012958,Li20243374,SHAO2020103166} for further details on adaptive diffusion-based and variational models.

To our knowledge, only a few papers address forward-backward parabolic problems with nonstandard growth conditions. Inspired by the idea in \cite{Oliver20031380,THANH20141403}, in this paper we study problem \eqref{P1}--\eqref{P3} using Sobolev regularization. Namely, we first investigate the existence and uniqueness of weak solutions to the homogeneous Neumann IBVP for the regularized equation
\begin{equation}
	u_t=\operatorname{div}\left(\left(\frac{1}{1+\left|\nabla u\right|^2}+\delta\left|\nabla u\right|^{p(x)-2}\right)\nabla u\right) + \varepsilon \Delta u_t , \notag
\end{equation}
where $\varepsilon\in (0,1)$. We then analyze the limit as $\varepsilon \rightarrow 0^{+}$ of the family $\{u^{\varepsilon}\}_{\varepsilon}$ of weak solutions to the above regularized problems and establish the existence of Young measure solutions to the original problem \eqref{P1}--\eqref{P3}; see Definition \ref{definition of P1--P3}. The uniqueness of such a solution is left as an open problem. We emphasize here that, although the treatment in \cite{Yin2003659} (corresponding to equation \eqref{P1} with constant exponent $1 < p < 2$; see also \cite{Demoulini1996376}) involves Young measures, the spirit of that approach, which does not rely on a regularization scheme motivated by physical considerations, differs from the present one. Indeed, the approximation procedure in \cite{Tong2025180,Yin2003659} combines the explicit time discretization for the evolution equation with variational methods used to solve the stationary problem obtained at each time step. The difficulties arising from the nonconvexity of the potential are addressed via relaxation theory (see \cite{Kinderlehrer1992001} for details), whereas the feasibility of our approach is based on the observation that the discretized counterpart of the damping term $\Delta u_t$ provides the convexity of the stationary problem even when the corresponding potential is nonconvex (see Lemma \ref{lem 1}). However, the existence of such solutions for problem \eqref{P1}--\eqref{P3} cannot be established following the method in \cite{Yin2003659}, and, consequently, our work also serves as a theoretical complement to the studies in \cite{GUIDOTTI20123226,Oliver20031380,Yin2003659}. Our main result is stated as follows:
\begin{thm}[Existence of Young measure solutions]\label{main thm}
	Let $f\in W^{1,p(\cdot)}(\Omega)$. Then problem \eqref{P1}--\eqref{P3} admits at least one Young measure solution.
\end{thm}

The rest of this paper is organized as follows. Section \ref{Preliminaries} recalls some preliminaries on variable exponent spaces and Young measures, and introduces a family of regularized problems approximating problem \eqref{P1}--\eqref{P3}, together with the definition of Young measure solutions. Section \ref{Weak solutions of the regularized problem} establishes the existence and uniqueness of Young measure solutions to the regularized problem. Finally, Section \ref{Proof of Theorem 1} presents the proof of Theorem \ref{main thm}.

\section{Preliminaries}\label{Preliminaries}

In this section, we present some definitions and results concerning Lebesgue and Sobolev spaces with variable exponents, as well as Young measures. We also define a Young measure solution to problem \eqref{P1}--\eqref{P3} and establish a family of corresponding regularized problems.

\subsection{Variable exponent spaces}

We first introduce the notations for variable exponent spaces. Let $\Omega$ be an open subset of $\mathbb{R}^N$. We denote by $\mathcal{P}(\Omega)$ the set of all Lebesgue measurable functions $p:\Omega\rightarrow [1,\infty]$. Functions $p\in \mathcal{P}(\Omega)$ are called variable exponents on $\Omega$. If $p^{+}<\infty$, we say that $p$ is a bounded variable exponent. For any $p\in \mathcal{P}(\Omega)$, we define $p^{\prime}\in \mathcal{P}(\Omega)$ by $\frac{1}{p(y)}+\frac{1}{p^{\prime}(y)}=1$, where $\frac{1}{\infty}:=0$. The function $p^{\prime}$ is called the dual variable exponent of $p$. 

Let $p\in\mathcal{P}(\Omega)$. We define the semimodular
\begin{equation}
	\varrho_{L^{p(\cdot)}(\Omega)}(u):=\int_{\Omega} \left|u(x)\right|^{p(x)} d x . \notag
\end{equation}
The variable exponent Lebesgue space $L^{p(\cdot)}(\Omega)$ is given by
\begin{equation}
	L^{p(\cdot)}(\Omega)=\left\{u:\Omega\rightarrow\mathbb{R} \text{ measurable}, \varrho_{L^{p(\cdot)}(\Omega)}(\lambda u)<\infty \text { for some } \lambda>0\right\} , \notag
\end{equation}
equipped with the Luxemburg norm
\begin{equation}
	\left\|u\right\|_{L^{p(\cdot)}(\Omega)}=\inf \left\{\lambda>0: \varrho_{L^{p(\cdot)}(\Omega)}\left(\frac{u}{\lambda}\right) \leq 1\right\} . \notag
\end{equation}
The function $u \in L^{p(\cdot)}(\Omega)$ belongs to the space $W^{k, p(\cdot)}(\Omega)$, where $k \in \mathbb{N}$ and $p \in \mathcal{P}(\Omega)$, if its weak partial derivatives $\partial_\alpha u$ with $\left|\alpha\right| \leq k$ exist and lie in $L^{p(\cdot)}(\Omega)$. We define the Sobolev semimodular on $W^{k,p(\cdot)}(\Omega)$ by
\begin{equation}
	\varrho_{W^{k, p(\cdot)}(\Omega)}(u):=\sum_{0 \leq\left|\alpha\right| \leq k} \varrho_{L^{p(\cdot)}(\Omega)}(\partial_\alpha u) , \notag
\end{equation}
which induces the norm
\begin{equation}
	\left\|u\right\|_{W^{k, p(\cdot)}(\Omega)}=\inf \left\{\lambda>0: \varrho_{W^{k, p(\cdot)}(\Omega)}\left(\frac{u}{\lambda}\right) \leq 1\right\} . \notag
\end{equation}
Furthermore, the variable exponent Sobolev space $W^{k, p(\cdot)}(\Omega)$ admits the equivalent norm
$$
\sum_{m=0}^k \left\|\nabla^m u\right\|_{L^{p(\cdot)}(\Omega)} = \sum_{m=0}^k\left\| \left|\nabla^m u\right| \right\|_{L^{p(\cdot)}(\Omega)} .
$$

We recall the log-H\"{o}lder continuity condition, which is widely used in the study of variable exponent spaces.
\begin{defn}
	A function $\beta: \Omega \rightarrow \mathbb{R}$ is called locally $\log$-H\"{o}lder continuous on $\Omega$ if there exists $a_1>0$ such that
	$$
	\left|\beta(x)-\beta(y)\right| \leq \frac{a_1}{\log \left(e+1/\left|x-y\right|\right)}
	$$
	for all $x, y \in \Omega$. A function $\beta$ is said to satisfy the log-H\"{o}lder decay condition if there exist $\beta_{\infty} \in \mathbb{R}$ and a constant $a_2>0$ such that
	$$
	\left|\beta(x)-\beta_{\infty}\right| \leq \frac{a_2}{\log \left(e+\left|x\right|\right)}
	$$
	for all $x \in \Omega$. We say that $\beta$ is globally $\log$-H\"{o}lder continuous in $\Omega$ if it is locally log-H\"{o}lder continuous and satisfies the log-H\"{o}lder decay condition. The constants $a_1$ and $a_2$ are referred to as the local log-H\"{o}lder constant and the log-H\"{o}lder decay constant, respectively. The maximum $\max \{a_1, a_2 \}$ is called the log-H\"{o}lder constant of $\beta$.
\end{defn}

We introduce the following class of variable exponents
$$
\mathcal{P}^{\log }(\Omega):=\left\{p \in \mathcal{P}(\Omega): \frac{1}{p} \text { is globally log-H\"{o}lder continuous}\right\} .
$$
The $\log$-H\"{o}lder constant of $\frac{1}{p}$ is denoted by $c_{\log}(p)$ or $c_{\log}$. Moreover, if $p \in \mathcal{P}(\Omega)$ with $p^{+}<\infty$, then $p \in \mathcal{P}^{\log}(\Omega)$ if and only if $p$ is globally $\log$-H\"{o}lder continuous.
\begin{defn}
	A bounded domain $\Omega \subset \mathbb{R}^N$ is called an $\alpha$-John domain, with $\alpha>0$, if there exists $x_0 \in \Omega$ such that each point in $\Omega$ can be joined to $x_0$ by a rectifiable path $\gamma$, parametrized by arc-length, satisfying
	$$
	B\left(\gamma(t), \frac{1}{\alpha} t\right) \subset \Omega
	$$
	for all $t \in[0, \ell(\gamma)]$, where $\ell(\gamma)$ is the length of $\gamma$. The ball $B\left(x_0, \frac{1}{2 \alpha} \operatorname{diam}(\Omega)\right)$ is called the John ball. If the value of $\alpha$ is irrelevant, we refer to $\Omega$ as a John domain.
	
	There are many examples of John domains. Every ball is clearly a 1-John domain, and every bounded domain satisfying the uniform interior cone condition is a John domain. In particular, any bounded Lipschitz domain is a John domain, where a domain with $C^{0,1}$-boundary is referred to as a Lipschitz domain; see also \cite{Diening2011,Martio1979383}.
\end{defn}

For later use, we state the following important properties of the spaces $L^{p(\cdot)}(\Omega)$ and $W^{1,p(\cdot)}(\Omega)$. Detailed discussions can be found in \cite{Diening2011,Kovavcik1991592,Mihuailescu20072929}.
\begin{lem}\cite{Diening2011}\label{norm vs semi norm}
	Let $p \in \mathcal{P}(\Omega)$ with $p^{-}<\infty$. If $\varrho_{L^{p(\cdot)}(\Omega)}(u)>0$ or $p^{+}<\infty$, then
	\begin{equation}\label{seminorm norm inequalities}
		\min \left\{\left\|u\right\|_{L^{p(\cdot)}(\Omega)}^{p^{-}}, \left\|u\right\|_{L^{p(\cdot)}(\Omega)}^{p^{+}}\right\} \leq \varrho_{L^{p(\cdot)}(\Omega)}(u) \leq \max \left\{\left\|u\right\|_{L^{p(\cdot)}(\Omega)}^{p^{-}}, \left\|u\right\|_{L^{p(\cdot)}(\Omega)}^{p^{+}}\right\} . \notag
	\end{equation}
\end{lem}
\begin{lem}\cite{Diening2011}\label{basic sobolev general trace}
	Let $p \in \mathcal{P}(\Omega)$. Then $L^{p(\cdot)}(\Omega)$ and $W^{k, p(\cdot)}(\Omega)$ are Banach spaces, separable if $p$ is bounded, and reflexive and uniformly convex whenever $1<p^{-} \leq p^{+}<\infty$.
\end{lem}
\begin{lem}[H\"{o}lder's inequality]\cite{Diening2011}\label{holder}
	Let $p, q, s \in \mathcal{P}(\Omega)$ be such that $\frac{1}{s(y)}=\frac{1}{p(y)}+\frac{1}{q(y)}$ for a.e. $y \in \Omega$. Then, for all $u \in L^{p(\cdot)}(\Omega)$ and $v \in L^{q(\cdot)}(\Omega)$,
	\begin{equation}
		\left\| uv \right\|_{L^{s(\cdot)}(\Omega)}\leq 2 \left\| u \right\|_{L^{p(\cdot)}(\Omega)}\left\| v \right\|_{L^{q(\cdot)}(\Omega)} , \notag
	\end{equation}
	where in the case $s=p=q=\infty$, we adopt the convention $\frac{s}{p}=\frac{s}{q}=1$.
\end{lem}
\begin{lem}\cite{Diening2011}\label{density sobolev general trace}
	Let $\Omega$ be a Lipschitz domain, and let $p \in \mathcal{P}^{\log }(\Omega)$ satisfy $1 \leq p^{-} \leq p^{+}<\infty$. Then $C^{\infty}(\bar{\Omega})$ is dense in $W^{k, p(\cdot)}(\Omega)$, $k \in \mathbb{N}$.
\end{lem}
\begin{lem}[Poincar\'{e} inequality] \cite{Diening2011} \label{Poincare inequality general trace}
	Let $p \in \mathcal{P}^{\log}(\Omega)$. If $\Omega$ is a bounded $\alpha$-John domain, then
	\begin{equation}
		\left\|u-\frac{1}{\left|\Omega\right|}\int_{\Omega} u dx\right\|_{L^{p(\cdot)}(\Omega)} \leq c \operatorname{diam}(\Omega)\left\|\nabla u\right\|_{L^{p(\cdot)}(\Omega)} , \notag
	\end{equation}
	for $u \in W^{1, p(\cdot)}(\Omega)$. The constant $c$ depends only on the dimension $N$, $\alpha$ and $c_{\log}(p)$.
\end{lem}

\begin{lem}\label{weak convergence lemma}
	Let $p\in \mathcal{P}(\Omega)$ with $1 < p^{-} \leq p^{+}<\infty$. Suppose $\{f^m\}_{m=1}^{\infty} \subset L^{\infty}(0,T;L^{p(\cdot)}(\Omega))$ and $f^m \rightharpoonup 0$ weakly $\ast$ in $L^{\infty}(0,T;L^{p(\cdot)}(\Omega))$ as $m \rightarrow \infty$. Let $\varPsi$ be a compact subset of $L^1(0,T;L^{p^{\prime}(\cdot)}(\Omega))$. Then
	\begin{equation}
		\sup _{\xi \in \varPsi}\left|\iint_{Q_T} f^m \xi dx dt\right| \rightarrow 0, \quad \text {as } m \rightarrow \infty.\notag
	\end{equation}
\end{lem}
\begin{proof}
	Here we may assume that $\left\|f^m\right\|_{L^{\infty}(0,T;L^{p(\cdot)}(\Omega))}\leq C$ for all $m\in \mathbb{N}$. Since $\varPsi\subset L^1(0,T;L^{p^{\prime}(\cdot)}(\Omega))$ is compact, it is totally bounded. Given any $r>0$, there exists a finite subset $\left\{\xi_1,\xi_2,\cdots,\xi_n\right\} \subset \varPsi$ such that, for every $\eta\in \varPsi$, one can find $i \in\{1,2, \cdots, n\}$ satisfying
	$$
	\left\| \eta - \xi_i \right\|_{L^1(0,T;L^{p^{\prime}(\cdot)}(\Omega))} < r . 
	$$
	By the weak $\ast$ convergence of $\{f^m\}_{m=1}^{\infty}$ (see \cite{Zeidler1990}), there exists $N_i\in\mathbb{N}$ such that, for all $m>N_i$, it holds that
	$$
	\left|\iint_{Q_T}f^m\xi_i dxdt\right| < r,\quad i=1,2,\cdots,n.
	$$
	Fix $N=\max_{1\leq i\leq n} N_i$. For any $\xi \in\varPsi$, there exists $j \in\{1,2, \cdots, n\}$ such that $\left\| \xi - \xi_j \right\|_{L^1(0,T;L^{p^{\prime}(\cdot)}(\Omega))} < r$. Consequently, for all $m>N$, it follows from Lemma \ref{holder} that
	\begin{align}
		\left|\iint_{Q_T} f^{m} \xi dx dt\right|\leq & \left| \iint_{Q_T} f^{m}\xi_{j}dxdt \right| + \left| \int_{0}^{T}\int_{\Omega} f^{m} \left(\xi - \xi_{j}\right) dxdt \right| \notag\\
		\leq & r + 2\int_{0}^{T} \left\| f^{m}(t) \right\|_{L^{p(\cdot)}(\Omega)} \left\| \xi(t) - \xi_{j}(t) \right\|_{L^{p^{\prime}(\cdot)}(\Omega)} dt \notag\\
		\leq & r + 2 \left\| f^{m} \right\|_{L^{\infty}(0,T;L^{p(\cdot)}(\Omega))} \left\| \xi - \xi_{j} \right\|_{L^{1}(0,T;L^{p^{\prime}(\cdot)}(\Omega))} \notag\\
		\leq & \left(1 + 2 C \right)r . \notag
	\end{align}
	This completes the proof.
\end{proof}

Adapting the argument of Theorem 3.20 in \cite{Dacorogna2008}, the following sufficient condition for weak lower semicontinuity on $W^{1,p(\cdot)}(\Omega)$ can be derived. The proof is omitted for brevity.
\begin{lem}\label{wlsc}
	Assume $p\in \mathcal{P}(\Omega)$ with $1 \leq p^{-} \leq p^{+}<\infty$, and let $f: \Omega \times \mathbb{R}^N \rightarrow$ $\mathbb{R} \cup\{+\infty\}$ be a Carath\'{e}odory function satisfying
	\begin{equation}
		f(x, \xi) \geq a(x) \cdot \xi + b(x)\notag
	\end{equation}
	for a.e. $x \in \Omega$ and for all $\xi \in \mathbb{R}^N$, where $a \in L^{p^{\prime}(\cdot)}(\Omega;\mathbb{R}^N)$ and $b \in L^1(\Omega)$. Let
	\begin{equation}
		J(\xi):=\int_{\Omega} f(x, \xi(x)) d x . \notag
	\end{equation}
	Assume that $\xi \rightarrow f(x, \xi)$ is convex and that
	\begin{equation}
		\xi_\nu \rightharpoonup \xi \text { weakly in } L^{p(\cdot)}(\Omega) . \notag
	\end{equation}
	Then
	\begin{equation}
		\liminf _{\nu \rightarrow \infty} J\left(\xi_\nu\right) \geq J(\xi) . \notag
	\end{equation}
\end{lem}

\subsection{Young measures}
For the readers' convenience, we briefly recall some fundamental notions concerning Young measures. Denote by $C_0(\mathbb{R}^N)$ the closure of continuous functions in $\mathbb{R}^N$ with compact support. The dual of $C_0(\mathbb{R}^N)$ can be identified with the space $\mathcal{M}(\mathbb{R}^N)$ of signed Radon measures with finite mass, via the pairing
$$
\langle\mu, f\rangle=\int_{\mathbb{R}^N} f(\xi) d \mu(\xi),\quad f\in C_0(\mathbb{R}^N),~\mu\in\mathcal{M}(\mathbb{R}^N).
$$
A map $\nu:\Omega\rightarrow\mathcal{M}(\mathbb{R}^N)$ is called weakly $\ast$ measurable if the functions $x\mapsto\langle\nu_x,f\rangle$ are measurable for all $f\in C_0(\mathbb{R}^N)$, where $\nu_x=\nu(x)$. We state the following fundamental theorem on Young measures, as given by J. M. Ball.
\begin{lem}\cite{ball2005version}\label{Fundamental theorem on Young measures}
	Let $\Omega \subset \mathbb{R}^N$ be Lebesgue measurable, let $K \subset \mathbb{R}^N$ be closed, and let $u_j: \Omega \rightarrow \mathbb{R}^N, j \in \mathbb{N}$, be a sequence of Lebesgue measurable functions satisfying $u_j \rightarrow K$ in measure as $j \rightarrow \infty$, i.e., given any open neighbourhood $U$ of $K$ in $\mathbb{R}^N$
	$$
	\lim _{j \rightarrow \infty}\left|\left\{x \in \Omega: u_j(x) \notin U\right\}\right|=0 .
	$$
	Then there exists a subsequence $\{u_{j_k}\}_{k=1}^{\infty}$ of $\{u_j\}_{j=1}^{\infty}$ and a weakly $\ast$ measurable map $\nu:\Omega\rightarrow\mathcal{M}(\mathbb{R}^N)$ such that
	\begin{itemize}
		\item [\emph{(i)}] $\nu(x) \geq 0$ and  $\|\nu(x)\|_{\mathcal{M}\left(\mathbb{R}^N\right)}=\int_{\mathbb{R}^N} d \nu_x \leq 1$ for a.e. $x \in \Omega$;
		\item [\emph{(ii)}] $\operatorname{spt} \nu_x \subset K$ for a.e. $x \in \Omega$;
		\item [\emph{(iii)}] For all $\psi\in C_0(\mathbb{R}^N)$, $\psi(u_{j_k})\rightharpoonup\langle\nu,\psi\rangle$ weakly $\ast$ in $L^{\infty}(\Omega)$.
	\end{itemize}
	Suppose further that $\{u_{j_k}\}_{k=1}^{\infty}$ satisfies the boundedness condition
	\begin{equation}\label{boundedness condition}
		\forall R>0:~\lim _{L \rightarrow \infty} \sup _{k \in \mathbb{N}} \left|\{x \in\Omega\cap B_R:|u_{j_k}(x)| \geq L\}\right|=0,
	\end{equation}
	where $B_R = B_R(0)$. Then
	\begin{equation}\label{norm nu}
		\left\|\nu_x\right\|_{\mathcal{M}(\mathbb{R}^N)}=1 \quad \text { for a.e. } x \in \Omega,
	\end{equation}
	and the following condition holds
	\begin{equation}\label{weak representation}
		\left\{\begin{array}{l}
			\text{For any measurable } D \subset \Omega \text { and any continuous function } f: \mathbb{R}^N \rightarrow \mathbb{R} \text { such } \\
			\text {that }\left\{f(u_{j_k})\right\} \text { is sequentially weakly relatively compact in } L^1(D), \text { we have } \\
			f(u_{j_k}) \rightharpoonup\left\langle\nu_x, f\right\rangle \text { in } L^1(D) \text {. }
		\end{array}\right.
	\end{equation}
\end{lem}
\begin{defn}
	The map $\nu: \Omega \rightarrow \mathcal{M}(\mathbb{R}^N)$ in Lemma \ref{Fundamental theorem on Young measures} is called the Young measure in $\mathbb{R}^N$ generated by the sequence $\{u_{j_k}\}_{k=1}^{\infty}$.
\end{defn}

Hungerb\"{u}hler subsequently established the following consequence, while refined versions of Lemma \ref{Fundamental theorem on Young measures} can be found in \cite{Kinderlehrer1994059}.
\begin{lem}\cite{Hungerbuhler1997}\label{equivalent lemma}
	Let $\Omega$, $u_j$, and $\nu_x$ be as in Lemma \ref{Fundamental theorem on Young measures}. Then \eqref{boundedness condition}, \eqref{norm nu} and \eqref{weak representation} are equivalent.
\end{lem}
\begin{rem}\label{rmk 1}
	In \cite{ball2005version}, it is also shown that, under hypothesis \eqref{boundedness condition}, for any measurable set $D \subset \Omega$,
	$$
	f(\cdot, u_{j_k}) \rightharpoonup\langle\nu_x, f(x, \cdot)\rangle\text { weakly in } L^1(D)
	$$
	for every Carath\'{e}odory function $f: D \times \mathbb{R}^N \rightarrow \mathbb{R}$ such that $\left\{f(\cdot, u_{j_k})\right\}_{k=1}^{\infty}$ is sequentially weakly relatively compact in $L^1(D)$. This statement is thus also equivalent to \eqref{boundedness condition}, \eqref{norm nu}, and \eqref{weak representation}.
\end{rem}

\subsection{Definition of Young measure solutions}
We define a Young measure solution to problem \eqref{P1}--\eqref{P3}. For clarity, let $\rho\in W^{1,p(\cdot)}(\Omega)$ be fixed but arbitrary, and define
\begin{equation}
	W_{\rho}^{1,p(\cdot)}(\Omega) = \left\{u\in W^{1,p(\cdot)}(\Omega):\int_{\Omega} u dx = \int_{\Omega} \rho dx\right\} \notag
\end{equation}
as the affine subspace of $W^{1,p(\cdot)}(\Omega)$. Recall that the variable exponent $p$ satisfies condition \eqref{p}, and denote by
\begin{equation}
	\vec{q}(x,\xi)=\frac{\xi}{1+\left|\xi\right|^2}+\delta\left|\xi\right|^{p(x)-2}\xi,\quad x\in\Omega,~\xi\in\mathbb{R}^N, \notag
\end{equation}
the heat flux in equation \eqref{P1}, where $\vec{q}(x,\xi) = \nabla_2 \varphi(x,\xi)$ and $\nabla_2$ denotes the derivative with respect to the second variable. The potential $\varphi:\Omega\times\mathbb{R}^N\rightarrow\mathbb{R}$ is a Carath\'{e}odory function, nonconvex with respect to $\xi$ for sufficiently small $\delta$. It is noted that $\varphi(x,\cdot)\in C^2(\mathbb{R}^N)$ for all $x\in\Omega$, and that $\vec{q}$ and $\varphi$ satisfy the following growth conditions:
\begin{equation}\label{structure condition 1}
	\left|\vec{q}(x,\xi)\right|\leq \lambda_2 \left|\xi\right|^{p(x)-1}+1,\quad\max\left\{\lambda_1 \left|\xi\right|^{p(x)}-1,0\right\} \leq\varphi(x,\xi) \leq \lambda_2 \left|\xi\right|^{p(x)}+1,
\end{equation}
for all $x\in\Omega$, $\xi\in\mathbb{R}^N$, where $0<\lambda_1\leq \lambda_2$ are constants. Furthermore, there exists a constant $K = K(N) > 0$ such that
\begin{equation}\label{A-B condition}
	\left(\vec{q}(x,\xi_1)-\vec{q}(x,\xi_2)\right)\cdot\left(\xi_1-\xi_2\right)\geq -K\left|\xi_1-\xi_2\right|^2,
\end{equation}
for all $x\in \Omega$, $\xi_1,\xi_2\in\mathbb{R}^N$.

\begin{defn}\label{definition of P1--P3}
	A pair $(u,\nu)$ is called a Young measure (or measure-valued) solution of problem \eqref{P1}--\eqref{P3} if $u\in L^{\infty}(0,T;W_f^{1,p(\cdot)}(\Omega))\cap H^1(0,T;L^2(\Omega))$ and $\nu = (\nu_{x,t})_{(x,t)\in Q_T}$ is a parametrized family of probability measures on $\mathbb{R}^N$, such that $u (0) = f$ a.e. in $\Omega$, and
	\begin{equation}
		\nabla u(x, t) = \langle\nu_{x, t},\mathrm{id}\rangle = \int_{\mathbb{R}^N} \xi d\nu_{x,t}(\xi) \quad \text{a.e. in }Q_T,\label{gradient}
	\end{equation}
	\begin{equation}
		\iint_{Q_T}\left(u_t\zeta + \langle\nu_{x,t},\vec{q}(x,\cdot) \rangle\cdot \nabla\zeta\right)dxdt=0,
	\end{equation}
	for every $\zeta\in L^{1}(0,T;W^{1,p(\cdot)}(\Omega))\cap L^2(0,T;H^1(\Omega))$, where $\mathrm{id}$ denotes the unit mapping in $\mathbb{R}^N$.
\end{defn}

Our goal is to construct a Young measure solution to problem \eqref{P1}--\eqref{P3}. Following the approach in \cite{Oliver20031380,THANH20141403}, we first consider the homogeneous Neumann IBVP for the Sobolev regularization of equation \eqref{P1}:
\begin{numcases}{}
	\displaystyle u^{\varepsilon}_t=\operatorname{div}\vec{q}(x,\nabla u^{\varepsilon})+\varepsilon\Delta u^{\varepsilon}_t, & $(x,t)\in Q_T$,\label{P4}\\
	\displaystyle \left(\vec{q}\left(x,\nabla u^{\varepsilon}\right)+\varepsilon\nabla u_t^{\varepsilon}\right)\cdot\vec{n}=0, & $(x,t)\in \partial\Omega\times (0,T)$,\label{P5} \\
	\displaystyle u^{\varepsilon}(x,0)=f, & $x\in\Omega$,\label{P6}
\end{numcases}
where $0<\varepsilon<1$. We refer to the problem \eqref{P4}--\eqref{P6} as the regularized problem corresponding to the original problem \eqref{P1}--\eqref{P3}.

\section{Weak solutions of the regularized problem \eqref{P4}--\eqref{P6}}\label{Weak solutions of the regularized problem}
In this section, we study the well-posedness of problem \eqref{P4}--\eqref{P6} and establish the following theorem by employing Rothe's method and variational principles.

\begin{thm}\label{thm 1}
	Let $f\in W^{1,p(\cdot)}(\Omega)$. Then the regularized problem \eqref{P4}--\eqref{P6} admits a unique weak solution $u^{\varepsilon}\in L^{\infty}(0,T;W_f^{1,p(\cdot)}(\Omega))\cap H^{1}(0,T;H^1(\Omega))$, satisfying $u^{\varepsilon}(0) = f$ a.e. in $\Omega$, and
	\begin{equation}\label{weak solution}
		\iint_{Q_T}\left(\frac{\partial u^{\varepsilon}}{\partial t}\zeta + \vec{q}(x,\nabla u^{\varepsilon})\cdot \nabla\zeta + \varepsilon \nabla \frac{\partial u^{\varepsilon}}{\partial t}\cdot\nabla\zeta \right)dxdt=0 ,
	\end{equation}
	for every $\zeta\in L^1(0,T;W^{1,p(\cdot)}(\Omega))\cap L^2(0,T;H^1(\Omega))$. Furthermore, for a.e. $t\in [0,T]$, the following inequality holds:
	\begin{equation}\label{solution estimate}
		\int_{\Omega}\varphi(x,\nabla u^{\varepsilon}(t))dx+\left\| \frac{\partial u^{\varepsilon}}{\partial \tau} \right\|^2_{L^2(Q_t)} + \left\| \nabla \sqrt{\varepsilon} \frac{\partial u^{\varepsilon}}{\partial \tau} \right\|^2_{L^2(Q_t)} \leq \int_{\Omega}\varphi(x,\nabla f) dx ,
	\end{equation}
	where $Q_t = \Omega \times (0,t)$.
	% where $M>0$ depends only on $N$, $\Omega$, $C$, $T$, $c$, $p$ $c_{\log}(p)$, and $f$, but is independent of $\varepsilon$.
\end{thm}
\begin{rem}
	By substituting the constant test function $\zeta = 1_{\Omega}$ into the weak formulation \eqref{solution estimate}, we obtain
	\begin{equation}
		\frac{d}{d t} \int_{\Omega} u^{\varepsilon} dx = 0,\notag
	\end{equation}
	which implies that the mean of any solution is conserved, i.e.,
	\begin{equation}
		\int_{\Omega} u^{\varepsilon}(x,t) dx = \int_{\Omega} u^{\varepsilon} (x,0) dx =  \int_{\Omega} f(x) dx . \notag
	\end{equation}
	Therefore, it is natural to require that the solution to problem \eqref{P4}--\eqref{P6} satisfies $u^{\varepsilon}(\cdot,t)\in W_f^{1,p(\cdot)}(\Omega)$ for a.e. $t\in [0,T]$.
\end{rem}

For simplicity, we omit the superscript $\varepsilon$ in the notation $u^{\varepsilon}$ and write $u$ instead. Let $m$ be a positive integer. For $u\in W^{1,p(\cdot)}(\Omega)$, we define
\begin{equation}
	J^{m,j}(u)=\int_{\Omega}\left( \varphi(x,\nabla u) + \frac{m}{2T}\left|u-u_m^{j-1}\right|^2 + \varepsilon\frac{m}{2T}\left|\nabla u - \nabla u_m^{j-1}\right|^2 \right) dx , \notag
\end{equation}
where $u_m^0=f$, $u_m^{j-1}\in W_f^{1,p(\cdot)}(\Omega)$, $j=1,2,\cdots,m$. The energy functional above is lower bounded due to the growth condition \eqref{structure condition 1}. It is readily seen that for every $v\in W^{1,p(\cdot)}(\Omega)$,
\begin{equation}
	J^{m,1}\left(v - \frac{1}{\left|\Omega\right|}\int_{\Omega} v dx + \frac{1}{\left|\Omega\right|}\int_{\Omega} f dx \right) \leq J^{m,1}(v) . \notag
\end{equation}
Without loss of generality, minimizers and minimizing sequences can be assumed to lie in $W_f^{1,p(\cdot)}(\Omega)$. First, we prove the following lemma.
\begin{lem}\label{lem 1}
	If $m$ is sufficiently large and $u_m^{j-1}\in W_f^{1,p(\cdot)}(\Omega)$, then $J^{m,j}$ admits a minimizer $u_m^j\in W_f^{1,p(\cdot)}(\Omega)$ for $j=1,2,\cdots,m$.
\end{lem}
\begin{proof}
	It is worth mentioning that the term in $J^{m, j}$ involving the discretized counterpart of $\nabla u_t$ ensures the convexity of the integrand with respect to $\nabla u$, although $\varphi$ itself is not convex. Introducing the notation
	\begin{equation}
		\varphi_{\varepsilon}(x,\xi) = \varphi(x,\xi) + \varepsilon\frac{m}{2T}\left|\xi\right|^2\notag
	\end{equation}
	and rewriting
	\begin{equation}
		J^{m,j}(u)=\int_{\Omega}\left( \varphi_{\varepsilon}(x,\nabla u) + \frac{m}{2T}\left|u-u_m^{j-1}\right|^2 -\varepsilon\frac{m}{T} \nabla u\cdot \nabla u_m^{j-1} + \varepsilon\frac{m}{2T}\left|\nabla u_m^{j-1}\right|^2 \right) dx . \notag
	\end{equation}
	It is clear that $u\mapsto J^{m,j}(u)-\int_{\Omega}\varphi_{\varepsilon}(x,\nabla u)dx$ is weakly lower semicontinuous on $W^{1,p(\cdot)}(\Omega)$. Additionally, by \eqref{A-B condition},
	\begin{equation}\label{treatment convexity}
		\left(\nabla_2\varphi_{\varepsilon}(x,\xi_1) - \nabla_2\varphi_{\varepsilon}(x,\xi_2)\right)\cdot\left(\xi_1-\xi_2\right)\geq \left(-K+\varepsilon\frac{m}{T}\right)\left|\xi_1-\xi_2\right|^2,\quad x\in \Omega,~\xi_1,\xi_2\in\mathbb{R}^N,
	\end{equation}
	which shows that the mapping $\xi\mapsto \varphi_{\varepsilon}(x,\xi)$ is convex for every $x\in \Omega$ whenever $m\geq\frac{KT}{\varepsilon}$. Lemma \ref{wlsc} then yields the weak lower semicontinuity of $u\mapsto \int_{\Omega}\varphi_{\varepsilon}(x,\nabla u)dx$ on $W^{1,p(\cdot)}(\Omega)$. Let $\{u_m^{j,k}\}_{k=1}^{\infty}\subset W_{f}^{1,p(\cdot)}(\Omega)$ be a minimizing sequence of $J^{m,j}$. Since $J^{m,j}$ is coercive on $W_{f}^{1,p(\cdot)}(\Omega)$ by Lemma \ref{Poincare inequality general trace}, we may assume that $\{u_m^{j,k}\}_{k=1}^{\infty}$ is bounded in $W^{1,p(\cdot)}(\Omega)$. Applying Lemma \ref{basic sobolev general trace}, there exists $u_m^j \in W_f^{1, p(\cdot)}(\Omega)$ and a subsequence of $\{u_m^{j, k}\}_{k=1}^{\infty}$, still denoted by itself, such that $u_m^{j,k}\rightharpoonup u_m^j$ weakly in $W^{1,p(\cdot)}(\Omega)$. Then, we have
	\begin{equation}
		\inf_{u\in W^{1,p(\cdot)}(\Omega)} J^{m,j}(u)\leq J^{m,j}(u_m^{j})\leq \liminf_{k \rightarrow \infty} J^{m,j}(u_m^{j,k}) = \inf_{u\in W^{1,p(\cdot)}(\Omega)} J^{m,j}(u) . \notag
	\end{equation}
	Consequently, $u_m^j \in W_f^{1, p(\cdot)}(\Omega)$ is a (not necessarily unique) minimizer of $J^{m,j}$. This finishes the proof.
	%The assertion follows from the direct method of the calculus of variations provided $J^{m,j}$ is weakly lower semicontinuous (wlsc) on $W^{1,p(\cdot)}(\Omega)$. 
\end{proof}

Note that, for $j=1,2,\cdots,m$, the minimizers $u_m^j$ satisfy the following Euler-Lagrange equations for every $\phi\in C^{\infty}(\bar{\Omega})$,
\begin{equation}\label{equ umj}
	\int_{\Omega}\left(\frac{m}{T}\left(u_m^{j}-u_m^{j-1}\right)\phi + \vec{q}(x,\nabla u_m^{j}) \cdot \nabla\phi + \varepsilon\frac{m}{T}\left(\nabla u_m^{j} - \nabla u_m^{j-1}\right)\cdot\nabla\phi \right)dx=0.
\end{equation}
From \eqref{structure condition 1} and the facts that $\nabla u_m^j\in L^{p(\cdot)}(\Omega)$ and $\vec{q}(x,\nabla u_m^j)\in L^{p^{\prime}(\cdot)}(\Omega)$, together with Lemma \ref{density sobolev general trace}, it follows that $\phi \in W^{1, p(\cdot)}(\Omega)$ is admissible as a test function in the the Euler-Lagrange system \eqref{equ umj} for $u_m^j$. We now derive an a priori estimate for the discrete energy.

\begin{lem}[Discrete energy inequality]\label{discrete energy inequality lemma}
	Given $\gamma\in (0,\varepsilon)$, there exists a constant $m_0(\gamma)>0$, such that for all $m\geq m_0$, the minimizers $u_m^j$ of the functionals $J^{m,j}$ satisfy the estimate
	\begin{equation}
		\int_{\Omega}\varphi(x,\nabla u_m^j) dx + \sum_{i=1}^{j}\int_{\Omega}\frac{m}{T}\left|u_m^i-u_m^{i-1}\right|^2 dx + \sum_{i=1}^{j}\int_{\Omega}(\varepsilon-\gamma)\frac{m}{T}\left|\nabla u_m^i-\nabla u_m^{i-1}\right|^2 dx \leq \int_{\Omega}\varphi(x,\nabla f) dx, \label{discrete energy inequality}
	\end{equation}
	where $j=1,2,\cdots,m$.
\end{lem}
\begin{proof}
	Similar to the treatment in \eqref{treatment convexity}, it follows from \eqref{A-B condition} that, for $m\geq \frac{KT}{2\gamma}$, the mapping $\xi\mapsto\varphi(x,\xi)+\gamma\frac{m}{T}|\xi-\nabla u_m^{j-1}|^2$, $j=1,2,\cdots,m$, is convex; i.e., for all $x\in\Omega$ and $\xi_1,\xi_2\in\mathbb{R}^N$, we have
	\begin{align}
		\varphi(x,\xi_1)+&\gamma\frac{m}{T}\left|\xi_1-\nabla u_m^{j-1}(x)\right|^2-\left(\varphi(x,\xi_2)+\gamma\frac{m}{T}\left|\xi_2-\nabla u_m^{j-1}(x)\right|^2\right)\notag\\
		&\leq \left(\vec{q}(x,\xi_1)+2\gamma\frac{m}{T}\left(\xi_1-\nabla u_m^{j-1}(x)\right)\right)\cdot\left(\xi_1-\xi_2\right) . \label{convexity conclu}
	\end{align}
	Set $m_0 = \max\left\{\frac{KT}{2\gamma},\frac{KT}{\varepsilon}\right\}$. Then, by Lemma \ref{lem 1}, $J^{m,j}$ admits a minimizer $u_m^j\in W_f^{1,p(\cdot)}(\Omega)$, satisfying \eqref{equ umj} for every $m\geq m_0$ and $j=1,2,\cdots,m$. From \eqref{convexity conclu}, for $i=1,2,\cdots,j$, we derive
	\begin{align}
		& \int_{\Omega}\varphi(x,\nabla u_m^{i})dx - \int_{\Omega}\varphi(x,\nabla u_m^{i-1})dx \notag\\
		= & \int_{\Omega} \bigg[ \varphi(x,\nabla u_m^{i}) + \gamma\frac{m}{T} \left|\nabla u_m^{i} - \nabla u_m^{i-1} \right|^2 - \left(\varphi(x,\nabla u_m^{i-1}) + \gamma\frac{m}{T} \left|\nabla u_m^{i-1} - \nabla u_m^{i-1} \right|^2\right) \notag\\
		& - \gamma\frac{m}{T} \left|\nabla u_m^{i} - \nabla u_m^{i-1} \right|^2
		\bigg] dx \notag\\
		\leq & \int_{\Omega} \bigg[ \left(\vec{q}(x,\nabla u_m^{i})+ 2\gamma\frac{m}{T} \left(\nabla u_m^{i} - \nabla u_m^{i-1}\right)\right)\cdot\left(\nabla u_m^{i} - \nabla u_m^{i-1}\right)\notag\\ 
		& - \gamma\frac{m}{T} \left|\nabla u_m^{i} - \nabla u_m^{i-1} \right|^2 \bigg] dx \notag\\
		= & \int_{\Omega} \bigg[ \left(\vec{q}(x,\nabla u_m^{i})+ \varepsilon\frac{m}{T} \left(\nabla u_m^{i} - \nabla u_m^{i-1}\right)\right)\cdot\left(\nabla u_m^{i} - \nabla u_m^{i-1}\right) \notag\\
		& - (\varepsilon-\gamma)\frac{m}{T} \left|\nabla u_m^{i} - \nabla u_m^{i-1} \right|^2 \bigg] dx. \label{discrete energy change}
	\end{align}
	Before proceeding with the estimate, we consider \eqref{equ umj} with $\phi=u_m^{i}-u_m^{i-1}$, which yields
	\begin{equation}
		\int_{\Omega}\left(\frac{m}{T}\left|u_m^{i}-u_m^{i-1}\right|^2 + \left(\vec{q}(x,\nabla u_m^{i}) + \varepsilon\frac{m}{T}\left(\nabla u_m^{i} - \nabla u_m^{i-1}\right)\right)\cdot \left(\nabla u_m^{i}-\nabla  u_m^{i-1}\right) \right)dx=0.\notag
	\end{equation}
	We insert the above equality into \eqref{discrete energy change} to obtain
	\begin{equation}
		\int_{\Omega}\varphi(x,\nabla u_m^{i})dx - \int_{\Omega}\varphi(x,\nabla u_m^{i-1})dx
		\leq -\int_{\Omega} \frac{m}{T}\left|u_m^{i}-u_m^{i-1}\right|^2 dx - \int_{\Omega} (\varepsilon-\gamma)\frac{m}{T} \left|\nabla u_m^{i} - \nabla u_m^{i-1} \right|^2 dx. \notag
	\end{equation}
	Summing the above inequality from $i=1$ to $j$ leads to
	\begin{equation}
		\int_{\Omega}\varphi(x,\nabla u_m^j) dx - \int_{\Omega}\varphi(x,\nabla f) dx \leq - \sum_{i=1}^{j}\int_{\Omega}\frac{m}{T}\left|u_m^i-u_m^{i-1}\right|^2 dx - \sum_{i=1}^{j}\int_{\Omega}(\varepsilon-\gamma)\frac{m}{T}\left|\nabla u_m^i-\nabla u_m^{i-1}\right|^2 dx. \notag
	\end{equation}
	This completes the proof.
\end{proof}

\begin{cor}\label{cor1}
	For every $m\geq m_{0}(\gamma)$, the minimizers of the functionals $J^{m,j}$ further satisfy the following estimate
	\begin{equation}
		\sup_{1\leq j\leq m}\left\| u_m^j \right\|_{W^{1,p(\cdot)}(\Omega)} + \sum_{j=1}^{m} \frac{m}{T} \left\| u_m^j - u_m^{j-1} \right\|_{L^2(\Omega)}^2 + \sum_{j=1}^{m} \left(\varepsilon - \gamma\right) \frac{m}{T} \left\|\nabla u_m^j - \nabla u_m^{j-1} \right\|_{L^2(\Omega)}^2\leq M_1, \label{discrete energy corollary}
	\end{equation}
	where $M_1>0$ depends only on $N$, $\Omega$, $\lambda_1$, $\lambda_2$, $p^{+}$, $c_{\log}(p)$, $f$, and $\alpha$, the John constant of $\Omega$.
\end{cor}

The assertion follows directly from the growth condition \eqref{structure condition 1}, Lemma \ref{Poincare inequality general trace} and \ref{discrete energy inequality lemma}, and the assumption that $f\in W^{1,p(\cdot)}(\Omega)$.

Now we construct approximate solutions to the regularized problem \eqref{P4}--\eqref{P6}. Let $m\geq m_{0}$. For $j=1,2,\cdots,m$, let $\chi_m^j(t)$ denote the characteristic function of the interval $[(j-1)\frac{T}{m},j\frac{T}{m})$, and define $\lambda_m^j(t)$ as
\begin{equation}
	\lambda_m^j(t)=\left(\frac{m}{T}t-(j-1)\right)\chi_m^j(t),\quad 0\leq t\leq T.\notag
\end{equation}
Define the following functions on $Q_T$
\begin{align}
	&u_m(x,t)=\sum_{j=1}^{m}\chi_m^j(t)\left\{u_m^{j-1}(x)+\lambda_m^j(t)\left(u_m^j(x)-u_m^{j-1}(x)\right)\right\},\notag\\
	&\bar{u}_m(x,t)=\sum_{j=1}^{m}\chi_m^j(t)u_m^j(x).\notag
\end{align}
Then the Euler-Lagrange equation \eqref{equ umj} can be equivalently rewritten as
\begin{equation}\label{eqn um}
	\iint_{Q_T}\left(\frac{\partial u_m}{\partial t}\zeta + \vec{q}(x,\nabla \bar{u}_m) \cdot \nabla\zeta + \varepsilon \nabla \frac{\partial u_m}{\partial t}\cdot\nabla\zeta \right)dxdt=0,
\end{equation}
for every $\zeta\in L^1(0,T;W^{1,p(\cdot)}(\Omega))\cap L^2(0,T;H^1(\Omega))$. To establish the existence of weak solutions to problem \eqref{P4}--\eqref{P6}, we use Corollary \ref{cor1} to obtain uniform estimates for the approximate solutions.
\begin{lem}
	There exists a positive constant $M_2$, depending only on $M_1$, $\Omega$, $T$, $\lambda_2$, and $p$, such that the following inequalities hold:
	\begin{align}
		& \sup_{0\leq t\leq T} \left\| u_m(t) \right\|_{W^{1,p(\cdot)}(\Omega)} \leq M_2,~\left\| \frac{\partial u_m}{\partial t} \right\|_{L^{2}(Q_T)} \leq M_2,~\left\|\nabla\sqrt{\varepsilon - \gamma}\frac{\partial u_m}{\partial t} \right\|_{L^{2}(Q_T)} \leq M_2, \label{estimates for um} \\
		& \sup_{0\leq t\leq T} \left\| \bar{u}_m(t) \right\|_{W^{1,p(\cdot)}(\Omega)} \leq M_2,~\sup_{0\leq t\leq T}\left\| \vec{q}(x, \nabla \bar{u}_m) \right\|_{L^{p^{\prime}(\cdot)}(\Omega)}\leq M_2, \label{estimates for barum} \\
		& \left\| u_m - \bar{u}_m \right\|_{L^2(Q_T)}\leq\frac{1}{m}M_2. \label{estimate for um minus barum} 
	\end{align}
\end{lem}
\begin{proof}
	For any $t\in [0,T)$, there exists a positive integer $j=1,2,\cdots,m$ such that $(j-1)\frac{T}{m}\leq t<j\frac{T}{m}$. Given that $\lambda_m^j(t)\in [0,1)$, it follows from \eqref{discrete energy corollary} that
	\begin{equation}
		\left\| u_m(t) \right\|_{W^{1,p(\cdot)}(\Omega)} \leq \lambda_m^j(t) \left\| u_m^j \right\|_{W^{1,p(\cdot)}(\Omega)} + \left(1 - \lambda_m^j(t)\right) \left\| u_m^{j-1} \right\|_{W^{1,p(\cdot)}(\Omega)} \leq M_1 , \notag
	\end{equation}
	and, analogously, we obtain the corresponding estimate for $\bar{u}_m$. Furthermore, we derive from \eqref{structure condition 1} and Lemma \ref{norm vs semi norm} that
	\begin{align}
		\varrho_{L^{p^{\prime}(\cdot)}(\Omega)}\left(\vec{q}(x,\nabla \bar{u}_m(t))\right) = &\int_{\Omega}\left|\vec{q}(x,\nabla u_m^{j})\right|^{\frac{p(x)}{p(x)-1}} dx \notag\\
		\leq & \int_{\Omega} \left(\lambda_2 \left|\nabla u_m^j\right|^{p(x)-1} + 1 \right)^{\frac{p(x)}{p(x)-1}} dx \notag\\
		\leq & \int_{\Omega} 2^{\frac{1}{p(x)-1}} \left( \lambda_{2}^{\frac{p(x)}{p(x)-1}}\left|\nabla u_m^j\right|^{p(x)} + 1 \right) dx \notag\\
		\leq & 2^{\frac{1}{p^{-}-1}} (\lambda_2+1)^{\frac{p^{-}}{p^{-}-1}} \int_{\Omega} \left|\nabla u_m^j\right|^{p(x)} dx + 2^{\frac{1}{p^{-}-1}} \left|\Omega\right| \notag\\
		\leq & 2^{\frac{1}{p^{-}-1}} \left(\lambda_2+1\right)^{\frac{p^{-}}{p^{-}-1}} \max \left\{\left\|\nabla u_m^j\right\|_{L^{p(\cdot)}(\Omega)}^{p^{-}},\left\|\nabla u_m^j\right\|_{L^{p(\cdot)}(\Omega)}^{p^{+}}\right\} + 2^{\frac{1}{p^{-}-1}}\left|\Omega\right| \notag\\
		\leq & 2^{\frac{1}{p^{-}-1}} \left(\lambda_2+1\right)^{\frac{p^{-}}{p^{-}-1}} \left(M_1 + 1\right)^{p^{+}} + 2^{\frac{1}{p^{-}-1}}\left|\Omega\right| \notag
	\end{align}
	which yields \eqref{estimates for barum}. Using \eqref{discrete energy corollary}, we obtain
	\begin{align}
		\left\| \nabla \sqrt{\varepsilon - \gamma}\frac{\partial u_m}{\partial t} \right\|^2_{L^{2}(Q_T)} = & \sum_{j=1}^m \left(\varepsilon - \gamma\right) \int_{(j-1)\frac{T}{m}}^{j\frac{T}{m}} \int_{\Omega} \left| \frac{m}{T} \left( \nabla u_m^{j} - \nabla u_m^{j-1}\right) \right|^2 dx dt \notag\\
		= & \sum_{j=1}^m \left(\varepsilon - \gamma\right) \frac{T}{m} \int_{\Omega} \left| \frac{m}{T} \left( \nabla u_m^{j} - \nabla u_m^{j-1}\right) \right|^2 dx \notag\\
		= & \sum_{j=1}^m \left(\varepsilon-\gamma\right) \frac{m}{T} \int_{\Omega} \left| \nabla u_m^{j} - \nabla u_m^{j-1} \right|^2 dx \leq M_1 . \notag
	\end{align}
	Similarly,
	\begin{align}
		\left\|\frac{\partial u_m}{\partial t}\right\|_{L^2(Q_T)}^2 = \sum_{j=1}^{m} \frac{m}{T} \int_{\Omega} \left|u_m^j - u_m^{j-1}\right|^2 dx \leq M_1 . \notag
	\end{align}
	Hence, we arrive at inequalities \eqref{estimates for um}. To prove \eqref{estimate for um minus barum}, it suffices to show that
	\begin{align}
		\left\| u_m - \bar{u}_m \right\|_{L^2(Q_T)}^2 = & \sum_{j=1}^m \int_{(j-1)\frac{T}{m}}^{j\frac{T}{m}} \int_{\Omega} \left| \left(\frac{m}{T}t-j\right)\left(u_m^j - u_m^{j-1}\right)\right|^2 dxdt \notag\\
		\leq & \sum_{j=1}^m \int_{(j-1)\frac{T}{m}}^{j\frac{T}{m}} \int_{\Omega} \left| u_m^j - u_m^{j-1}\right|^2 dxdt \notag \\
		= & \left(\frac{T}{m}\right)^2 \sum_{j=1}^m \frac{m}{T} \int_{\Omega} \left| u_m^j - u_m^{j-1}\right|^2 dx. \notag
	\end{align}
	The proof is thus complete.
\end{proof}

By Lemma \ref{basic sobolev general trace}, upon extracting suitable (not relabeled) subsequences, we obtain the following convergences as $m\to\infty$:
\begin{equation}\label{convergence um}
	\begin{aligned}
		& u_m \rightharpoonup u \text{ weakly }{\ast}\text{ in } L^{\infty}(0,T;W^{1,p(\cdot)}(\Omega)) \\
		& u_m \rightharpoonup u \text{ weakly in } H^1(0,T;H^1(\Omega)) \\
		& u_m \rightarrow u \text{ strongly in } C([0,T];L^{2}(\Omega)) \\
		& \sqrt{\varepsilon-\gamma}\frac{\partial u_m}{\partial t}\rightharpoonup \sqrt{\varepsilon-\gamma}\frac{\partial u}{\partial t}\text{ weakly in } L^2(0,T;H^1(\Omega)),\\
		& \bar{u}_m \rightharpoonup \bar{u} \text{ weakly }{\ast}\text{ in } L^{\infty}(0,T;W^{1,p(\cdot)}(\Omega)), \\
		& \vec{q}(x,\nabla \bar{u}_m) \rightharpoonup \vec{\Lambda} \text{ weakly }{\ast}\text{ in } L^{\infty}(0,T;L^{p^{\prime}(\cdot)}(\Omega)) . 
	\end{aligned}
\end{equation}
Since $u_m-\bar{u}_m$ converges strongly to zero in $L^2(Q_T)$, we deduce that $u=\bar{u}$ a.e. in $Q_T$. Moreover, $u(t) \in W_{f}^{1, p(\cdot)}(\Omega)$ a.e. in $[0,T]$ follows from \eqref{convergence um}. Passing to the limit as  $m\rightarrow\infty$ in \eqref{eqn um}, we then obtain
\begin{equation}\label{eqn Lambda}
	\iint_{Q_T}\left(\frac{\partial u}{\partial t}\zeta + \vec{\Lambda} \cdot \nabla\zeta + \varepsilon \nabla \frac{\partial u}{\partial t}\cdot\nabla\zeta \right)dxdt=0,
\end{equation}
for every $\zeta\in L^1(0,T;W^{1,p(\cdot)}(\Omega))\cap L^2(0,T;H^1(\Omega))$. To determine the function $\vec{\Lambda}$, we further consider the stronger convergence of the sequence $\{\nabla \bar{u}_m\}_{m\geq m_0}$, as stated in the following lemma.
\begin{lem}\label{strong convergence gradient um}
	Let $m$ denote the index of the subsequence along which the convergences in \eqref{convergence um} hold. Then both $\nabla u_m$ and $\nabla \bar{u}_m$ converge strongly in $L^2(0,T;L^2(\Omega))$ as $m\to\infty$.
\end{lem}
\begin{proof}
	We first examine the relationship between the piecewise constant and piecewise linear interpolations of a function defined at discrete time steps. 
	
	For any $i=1,2,\cdots,m$ and $t\in [(i-1)\frac{T}{m},i\frac{T}{m})$, by construction of approximate solutions, we have
	\begin{equation}
		\nabla u_m(t) - \nabla u(t) = \nabla u_m^i - \nabla u(t) - \left(i-\frac{m}{T}t\right)\left(\nabla u_m^i - \nabla u_m^{i-1}\right) , \notag
	\end{equation}
	and thus
	\begin{align}
		\int_{\Omega} \left|\nabla u_m^i- \nabla u(t)\right|^2 dx = & \int_{\Omega} \left|\nabla u_m(t) - \nabla u(t) + \left(i-\frac{m}{T}t\right)\left(\nabla u_m^i - \nabla u_m^{i-1}\right)\right|^2  dx \notag\\
		\leq & 2 \int_{\Omega} \left|\nabla u_m(t) - \nabla u(t)\right|^2 dx + 2 \left(i-\frac{m}{T}t\right)^2 \int_{\Omega} \left|\nabla u_m^i - \nabla u_m^{i-1}\right|^2 dx . \label{minus formula}
	\end{align}
	Integrating the above inequality over $[(i-1)\frac{T}{m}, i\frac{T}{m})$ gives
	\begin{align}
		& \int_{(i-1)\frac{T}{m}}^{i\frac{T}{m}} \int_{\Omega} \left|\nabla u_m^i- \nabla u\right|^2 dx dt \notag\\ \leq & 2 \int_{(i-1)\frac{T}{m}}^{i\frac{T}{m}} \int_{\Omega} \left|\nabla u_m - \nabla u\right|^2 dx dt + 2 \int_{(i-1)\frac{T}{m}}^{i\frac{T}{m}} \left(i-\frac{m}{T}t\right)^2 dt \int_{\Omega} \left|\nabla u_m^i - \nabla u_m^{i-1}\right|^2 dx \notag\\
		= & 2 \int_{(i-1)\frac{T}{m}}^{i\frac{T}{m}} \int_{\Omega} \left|\nabla u_m - \nabla u\right|^2 dx dt + \frac{2T}{3m} \int_{\Omega} \left|\nabla u_m^i - \nabla u_m^{i-1}\right|^2 dx . \label{L2norm gradientum minus gradientu}
	\end{align}
	For any $t\in [0,T)$, there exists an index $l\in \{1,2,\cdots,m\}$ such that $(l-1)\frac{T}{m}\leq t<l\frac{T}{m}$. For the case $l\geq 2$, from \eqref{discrete energy corollary}, \eqref{L2norm gradientum minus gradientu}, and the H\"{o}lder inequality (see Lemma \ref{holder}), we derive
	\begin{align}
		\int_{0}^{t}\int_{\Omega}\left|\nabla \bar{u}_m - \nabla u \right|^2 dx d\tau = & \sum_{i=1}^{l-1}\int_{(i-1)\frac{T}{m}}^{i\frac{T}{m}}\int_{\Omega} \left|\nabla u_m^i - \nabla u \right|^2 dx dt + \int_{(l-1)\frac{T}{m}}^{t}\int_{\Omega} \left|\nabla \bar{u}_m - \nabla u \right|^2 dx d\tau \notag\\
		\leq & \sum_{i=1}^{l-1}\left[ 2 \int_{(i-1)\frac{T}{m}}^{i\frac{T}{m}} \int_{\Omega} \left|\nabla u_m - \nabla u\right|^2 dx dt + \frac{2T}{3m} \int_{\Omega} \left|\nabla u_m^i - \nabla u_m^{i-1}\right|^2 dx \right] \notag\\
		& + \int_{(l-1)\frac{T}{m}}^{t}\int_{\Omega} \left|\nabla \bar{u}_m - \nabla u \right|^2 dx d\tau \notag\\
		\leq & 2 \int_{0}^{t} \int_{\Omega} \left|\nabla u_m - \nabla u\right|^2 dx d\tau + \frac{2T}{3m} \sum_{i=1}^{l-1} \int_{\Omega} \left|\nabla u_m^i - \nabla u_m^{i-1}\right|^2 dx \notag\\
		& + 2\int_{(l-1)\frac{T}{m}}^{t} \left\|\nabla \bar{u}_m(\tau) \right\|_{L^2(\Omega)}^2 d\tau + 2\int_{(l-1)\frac{T}{m}}^{t} \left\|\nabla u(\tau) \right\|_{L^2(\Omega)}^2 d\tau \notag\\
		\leq & 2 \int_{0}^{t} \int_{\Omega} \left|\nabla u_m - \nabla u\right|^2 dx d\tau + \frac{2T^2}{3m^2\left(\varepsilon-\gamma\right)} \sum_{i=1}^{l-1} \left(\varepsilon-\gamma\right) \frac{m}{T} \int_{\Omega} \left|\nabla u_m^i - \nabla u_m^{i-1}\right|^2 dx \notag\\
		& + 2\int_{(l-1)\frac{T}{m}}^{t} \left(2\left(1+\left|\Omega\right|\right)\right)^2 \left\|\nabla u_m(\tau) \right\|_{L^{p(\cdot)}(\Omega)}^2 d\tau \notag\\
		& + 2\int_{(l-1)\frac{T}{m}}^{t} \left(2\left(1+\left|\Omega\right|\right)\right)^2 \left\|\nabla u(\tau) \right\|_{L^{p(\cdot)}(\Omega)}^2 d\tau \notag\\
		\leq & 2 \int_{0}^{t} \int_{\Omega} \left|\nabla u_m - \nabla u\right|^2 dx d\tau + \frac{2 T^2}{3 m^2\left(\varepsilon-\gamma\right)}C_1 + \frac{T}{m} C_1, \label{relationship gradient um barum}
	\end{align}
	where $C_1>0$ depends only on $\Omega$, $M_1$, and $M_2$. As for the case $l=1$, the proof details are omitted. Indeed, integrating both sides of \eqref{minus formula} with $i=1$ over $(0,t)$ yields a result analogous to \eqref{relationship gradient um barum}.

	We now analyze the convergence of $\{\nabla u_m\}_{m\geq m_0}$ and $\{\nabla \bar{u}_m\}_{m\geq m_0}$ based on the preceding discussion.	By choosing the test functions $\zeta=\bar{u}_m-u$ in \eqref{eqn um} and $u_m-u$ in \eqref{eqn Lambda}, respectively, and then subtracting the resulting identities, we have, for $t\in [0,T]$,
	\begin{align}
		\int_{0}^{t}&\int_{\Omega} \left( \vec{q}(x,\nabla \bar{u}_m)\cdot\left(\nabla \bar{u}_m - \nabla u\right) - \vec{\Lambda}\cdot\left(\nabla u_m - \nabla u \right)\right) dxd\tau \notag\\
		&+\varepsilon\int_{0}^{t}\int_{\Omega} \left(\nabla \frac{\partial u_m}{\partial \tau}\cdot\left(\nabla\bar{u}_m - \nabla u\right)-\nabla \frac{\partial u}{\partial \tau}\cdot\left(\nabla u_m - \nabla u\right)\right)dxd\tau \notag\\
		&+\int_{0}^{t}\int_{\Omega} \left(\frac{\partial u_m}{\partial \tau}\left(\bar{u}_m - u\right)-\frac{\partial u}{\partial \tau}\left(u_m - u\right)\right)dxd\tau = 0 . \notag
	\end{align}
	We denote the three terms on the left-hand side of the above equality by $T_1$, $T_2$, and $T_3$, and write $\theta(m)$ for quantities tending to zero as $m\rightarrow\infty$, uniformly with respect to $t$. By \eqref{A-B condition} and Lemma \ref{weak convergence lemma} (e.g., with $\varPsi = \{\chi_{Q_t}\vec{q}(x,\nabla u):t\in [0,T]\}$), we obtain
	\begin{align}
		T_1 = & \int_{0}^{t}\int_{\Omega}  \left(\vec{q}(x,\nabla \bar{u}_m) - \vec{q}(x,\nabla u)\right)\cdot\left(\nabla \bar{u}_m - \nabla u\right) dxd\tau - \int_{0}^{t}\int_{\Omega} \vec{\Lambda}\cdot\left(\nabla u_m - \nabla u \right) dxd\tau \notag\\
		& + \int_{0}^{t}\int_{\Omega}  \vec{q}(x,\nabla u) \cdot \left(\nabla \bar{u}_m - \nabla u\right) dxd\tau \notag\\
		\geq & -K\int_{0}^{t}\int_{\Omega}\left|\nabla \bar{u}_m - \nabla u\right|^2 dxd\tau - \sup_{0\leq t\leq T}\left|\iint_{Q_T} \left(\chi_{Q_t}\vec{\Lambda}\right)\cdot\left(\nabla u_m - \nabla u \right) dxd\tau\right| \notag\\
		& - \sup_{0\leq t\leq T}\left|\iint_{Q_T} \left(\chi_{Q_t}\vec{q}(x,\nabla u)\right) \cdot \left(\nabla \bar{u}_m - \nabla u\right) dxd\tau\right| \notag\\
		\geq & -K\int_{0}^{t}\int_{\Omega}\left|\nabla \bar{u}_m - \nabla u\right|^2 dxd\tau - \theta(m) , \notag
	\end{align}
	which, together with \eqref{relationship gradient um barum}, implies that
	\begin{equation}
		T_1 \geq -2K \int_{0}^{t} \int_{\Omega} \left|\nabla u_m - \nabla u\right|^2 dx d\tau - \frac{\theta(m)}{\varepsilon-\gamma} - \theta(m) . \notag
	\end{equation}
	Direct calculation shows that
	\begin{align}
		T_2 = & \varepsilon \int_{0}^{t}\int_{\Omega}\left( \left(\nabla\frac{\partial u_m}{\partial \tau}-\nabla\frac{\partial u}{\partial \tau}\right)\cdot\left(\nabla u_m - \nabla u\right) + \nabla\frac{\partial u_m}{\partial \tau}\cdot\left(\nabla\bar{u}_m - \nabla u_m\right) \right) dxd\tau \notag\\
		= & \varepsilon \int_{0}^{t}\int_{\Omega} \left(\nabla\frac{\partial u_m}{\partial \tau}-\nabla\frac{\partial u}{\partial \tau}\right)\cdot\left(\nabla u_m - \nabla u\right) dxd\tau + \varepsilon\int_{0}^{t}\int_{\Omega} \nabla\frac{\partial u_m}{\partial \tau}\cdot\left(\nabla\bar{u}_m - \nabla u_m\right) dxd\tau \notag\\
		= & \frac{\varepsilon}{2} \int_{\Omega}\left|\nabla u_m(t)-\nabla u(t)\right|^2 dx - \frac{\varepsilon}{2} \int_{\Omega}\left|\nabla u_m(0)-\nabla u(0)\right|^2 dx \notag\\
		& + \varepsilon\int_{0}^{t}\int_{\Omega} \nabla\frac{\partial u_m}{\partial \tau}\cdot\left(\nabla\bar{u}_m - \nabla u_m\right) dxd\tau . \label{estimate of T2}
	\end{align}	
	It remains to estimate the last two terms on the right-hand side of \eqref{estimate of T2}. On the one hand, for any $t \in [0,T)$, assume that $t\in[(j-1)\frac{T}{m},j\frac{T}{m})$ for some $j=1,2,\cdots,m$. Then
	\begin{align}
		& \varepsilon\left|\int_{0}^{t}\int_{\Omega} \nabla\frac{\partial u_m}{\partial \tau}\cdot\left(\nabla\bar{u}_m - \nabla u_m\right) dxd\tau\right|\notag\\
		\leq & \varepsilon\int_{0}^{j\frac{T}{m}}\int_{\Omega} \left|\nabla\frac{\partial u_m}{\partial t}\cdot\left(\nabla\bar{u}_m - \nabla u_m\right)\right| dxdt\notag\\
		= & \varepsilon \sum_{i=1}^{j}\int_{(i-1)\frac{T}{m}}^{i\frac{T}{m}}\int_{\Omega} \frac{m}{T}\left(i-\frac{m}{T}t\right)\left|\nabla u_m^i - \nabla u_m^{i-1}\right|^2 dxdt \notag\\
		= & \varepsilon \sum_{i=1}^{j}\int_{(i-1)\frac{T}{m}}^{i\frac{T}{m}}\frac{m}{T}\left(i-\frac{m}{T}t\right)dt \int_{\Omega}\left|\nabla u_m^i - \nabla u_m^{i-1}\right|^2 dx \notag\\
		= & \frac{T\varepsilon}{2m\left(\varepsilon-\gamma\right)} \sum_{i=1}^{j}\left(\varepsilon-\gamma\right) \frac{m}{T}\int_{\Omega}\left|\nabla u_m^{i} - \nabla u_m^{i-1}\right|^2 dx \notag\\
		\leq & \frac{TM_1\varepsilon}{2m\left(\varepsilon-\gamma\right)}. \label{RHS third}
	\end{align}
	On the other hand, we choose the test function $\zeta = \xi\eta$ in equation \eqref{eqn um}, where $\xi\in C_c^1(\Omega)$ and $\eta\in C^1([0,T])$ satisfies $\eta(T) = 0$ and $\eta(0)\neq 0$. Using integration by parts, we get
	\begin{equation}
		-\int_{\Omega} f \xi\eta(0) dx-\iint_{Q_T} u_m \xi\eta_t dxdt + \iint_{Q_T} \left(\vec{q}(x,\nabla \bar{u}_m) \cdot \nabla\xi \eta + \varepsilon \nabla \frac{\partial u_m}{\partial t}\cdot\nabla\xi\eta \right)dxdt=0.\notag
	\end{equation}
	Letting $m\rightarrow\infty$, we have from \eqref{convergence um} that
	\begin{equation}
		-\int_{\Omega} f \xi\eta(0) dx-\iint_{Q_T} u \xi\eta_t dxdt + \iint_{Q_T}\left(\vec{\Lambda} \cdot \nabla\xi\eta + \varepsilon \nabla \frac{\partial u}{\partial t}\cdot\nabla\xi\eta \right)dxdt=0.\notag
	\end{equation}
	Similarly, substituting the test function $\zeta$ into \eqref{eqn Lambda}, we obtain
	\begin{equation}
		\iint_{Q_T}\left(\frac{\partial u}{\partial t}\xi\eta + \vec{\Lambda} \cdot \nabla\xi\eta + \varepsilon \nabla \frac{\partial u}{\partial t}\cdot\nabla\xi\eta \right)dxdt=0 . \notag
	\end{equation}
	Then
	\begin{equation}
		-\int_{\Omega}u(0)\xi\eta(0) dx - \iint_{Q_T}u\xi\eta_t dxdt + \iint_{Q_T}\left(\vec{\Lambda} \cdot \nabla\xi\eta + \varepsilon \nabla \frac{\partial u}{\partial t}\cdot\nabla\xi\eta \right)dxdt=0,\notag
	\end{equation}
	which yields $u(0)=f$ a.e. in $\Omega$. In summary, combining \eqref{estimate of T2} and \eqref{RHS third}, we deduce
	\begin{equation}
		T_2 \geq \frac{\varepsilon}{2} \int_{\Omega}\left|\nabla u_m(t)-\nabla u(t)\right|^2 dx - \frac{\varepsilon}{\varepsilon-\gamma} \theta(m) . \notag
	\end{equation}
	For the estimate of $T_3$, we apply H\"{o}lder's inequality and derive
	\begin{align}
		\left| T_3 \right| = & \left| \int_{0}^{t}\int_{\Omega} \left( \left(\frac{\partial u_m}{\partial \tau} - \frac{\partial u}{\partial \tau}\right)\left(u_m - u\right) + \frac{\partial u_m}{\partial \tau}\left(\bar{u}_m - u_m\right)\right) dxd\tau \right| \notag\\
		\leq & \frac{1}{2}\int_{\Omega}\left|u_m(t) - u(t)\right|^2 dx + \left|\int_{0}^{t}\int_{\Omega} \frac{\partial u_m}{\partial \tau}\left(\bar{u}_m - u_m\right) dxd\tau\right| \notag \\
		\leq & \frac{1}{2}\sup_{0\leq t\leq T}\left\| u_m(t) - u(t) \right\|_{L^2(\Omega)}^2 + \left\| \frac{\partial u_m}{\partial t} \right\|_{L^{2}(Q_T)} \left\| \bar{u}_m - u_m \right\|_{L^{2}(Q_T)} \notag \\
		\leq & \frac{1}{2} \left\| u_m - u \right\|_{L^{\infty}(0,T;L^{2}(\Omega))}^2 + \frac{1}{m} M_2^2 . \notag
	\end{align}
	Hence, from \eqref{estimates for um}, \eqref{estimate for um minus barum}, and \eqref{convergence um}, it follows that $T_3 = \theta(m)$. Taking into account that $T_1 + T_2 + T_3 = 0$ and combining this with the above estimates, we obtain
	\begin{equation}
		\frac{\varepsilon}{2}\int_{\Omega}\left|\nabla u_m(t) - \nabla u(t)\right|^2 dx \leq 2K\int_{0}^{t}\int_{\Omega}\left|\nabla u_m - \nabla u\right|^2 dxd\tau + \frac{\theta(m)}{\varepsilon-\gamma}+\theta(m) , \notag
	\end{equation}
	where $\theta(m)$ may be taken non-negative without loss of generality, otherwise we consider its absolute value. As a result,
	\begin{equation}
		\partial_t\int_{0}^{t}\int_{\Omega}\left|\nabla u_m - \nabla u\right|^2 dxd\tau\leq \frac{4K}{\varepsilon}\int_{0}^{t}\int_{\Omega}\left|\nabla u_m - \nabla u\right|^2 dxd\tau + \frac{\theta(m)}{\varepsilon\left(\varepsilon-\gamma\right)} . \notag
	\end{equation}
	Now, applying Gronwall's inequality, we arrive at
	\begin{equation}
		\int_{0}^{T}\int_{\Omega}\left|\nabla u_m - \nabla u\right|^2 dxdt\leq \frac{\theta(m)}{4K\left(\varepsilon-\gamma\right)}\exp\left(\frac{4KT}{\varepsilon}\right), \notag
	\end{equation}
	which tends to zero as $m\to\infty$. Therefore,
	\begin{equation}
		\nabla u_m \rightarrow \nabla u \text{ strongly in } L^2(0,T;L^2(\Omega)). \notag
	\end{equation}
	From \eqref{relationship gradient um barum}, we further deduce that 
	\begin{equation}\label{strong convergen barum}
		\nabla \bar{u}_m \rightarrow \nabla u \text{ strongly in } L^2(0,T;L^2(\Omega)).
	\end{equation}
	The proof of Lemma \ref{strong convergence gradient um} is complete.
\end{proof}

Next, we proceed to establish the form of the vector-valued function $\vec{\Lambda}$.

\begin{lem}
	It holds that
	\begin{equation}
		\vec{\Lambda}(x,t) = \vec{q}\left(x,\nabla u(x,t)\right) \text{ a.e. in } Q_T.\label{Lambda equal}
	\end{equation} 
\end{lem}
\begin{proof}
	Observe that for any $t\in [0,T]$, we have
	\begin{align}
		\varrho_{L^{p^{\prime}(\cdot)}(\Omega)}\left(\left|\nabla \bar{u}_m(t)\right|^{p(x)-2}\nabla \bar{u}_m(t)\right) = & \int_{\Omega}\left| \left|\nabla \bar{u}_m(t)\right|^{p(x)-2}\nabla \bar{u}_m(t) \right|^{\frac{p(x)}{p(x)-1}} dx \notag \\
		= & \int_{\Omega} \left|\nabla \bar{u}_m(t) \right|^{p(x)}dx \notag\\
		\leq & \max \left\{\left\| \bar{u}_m(t)\right\|_{W^{1,p(\cdot)}(\Omega)}^{p^{-}}, \left\| \bar{u}_m(t)\right\|_{W^{1,p(\cdot)}(\Omega)}^{p^{+}}\right\} , \notag
	\end{align}
	which gives
	\begin{equation}
		\sup_{0\leq t\leq T} \left\| \left|\nabla \bar{u}_m(t)\right|^{p(x)-2}\nabla \bar{u}_m(t) \right\|_{L^{p^{\prime}(\cdot)}(\Omega} \leq \left(M_2 + 1\right)^{p^{+}-1} . \notag
	\end{equation}
	Therefore, there exists a vector-valued function $\vec{\sigma}(x,t)$ such that
	\begin{equation}
		\left|\nabla \bar{u}_m\right|^{p(x)-2}\nabla \bar{u}_m \rightharpoonup \vec{\sigma} \text{ weakly }{\ast}\text{ in } L^{\infty}(0,T;L^{p^{\prime}(\cdot)}(\Omega) . \notag
	\end{equation}
	Considering the global Lipschitz continuity of the mapping $\xi\mapsto \frac{\xi}{1+\left|\xi\right|^2}$, we deduce from \eqref{strong convergen barum} that
	\begin{equation}
		\frac{\nabla \bar{u}_m}{1+\left|\nabla \bar{u}_m\right|^2} \rightarrow \frac{\nabla u}{1+\left|\nabla u\right|^2} \text{ strongly in }L^2(Q_T) . \label{strong convergence of frac grabarum}
	\end{equation}
	Hence, to prove \eqref{Lambda equal}, it suffices to show that
	$$
	\vec{\sigma}(x,t) = \left|\nabla u(x,t)\right|^{p(x)-2} \nabla u(x,t) \text{ a.e. in } Q_T .
	$$
	The integral identity \eqref{eqn um} can be rewritten as
	\begin{equation}\label{Lambda equal um}
		\iint_{Q_T}\left(\frac{\partial u_m}{\partial t}\zeta + \frac{\nabla \bar{u}_m}{1+\left|\nabla \bar{u}_m\right|^2}\cdot \nabla\zeta+\delta\left|\nabla \bar{u}_m\right|^{p(x)-2}\nabla \bar{u}_m\cdot \nabla\zeta + \varepsilon \nabla \frac{\partial u_m}{\partial t}\cdot\nabla\zeta \right)dxdt=0 ,
	\end{equation}
	for every $\zeta\in L^1(0,T;W^{1,p(\cdot)}(\Omega))\cap L^2(0,T;H^1(\Omega))$. Letting $m\rightarrow\infty$ in \eqref{Lambda equal um} and using \eqref{strong convergence of frac grabarum}, we have
	\begin{equation}\label{Lambda equal um sigma}
		\iint_{Q_T}\left(\frac{\partial u}{\partial t}\zeta + \frac{\nabla u}{1+\left|\nabla u\right|^2}\cdot \nabla\zeta+\delta\vec{\sigma}\cdot \nabla\zeta + \varepsilon \nabla \frac{\partial u}{\partial t}\cdot\nabla\zeta \right)dxdt=0 .
	\end{equation}
	Taking $\zeta = u$, we obtain
	\begin{equation}
		\iint_{Q_T}\left(\frac{\partial u}{\partial t} u + \frac{\left|\nabla u\right|^2}{1+\left|\nabla u\right|^2}+\delta\vec{\sigma}\cdot \nabla u + \varepsilon \nabla \frac{\partial u}{\partial t}\cdot\nabla u \right)dxdt=0 . \notag
	\end{equation}
	Furthermore, substituting $\zeta = \bar{u}_m$ into \eqref{Lambda equal um} yields
	\begin{equation}\label{Lambda equal um px}
		\iint_{Q_T}\left(\frac{\partial u_m}{\partial t}\bar{u}_m + \frac{\left|\nabla \bar{u}_m\right|^2}{1+\left|\nabla \bar{u}_m\right|^2}+\delta\left|\nabla \bar{u}_m\right|^{p(x)} + \varepsilon \nabla \frac{\partial u_m}{\partial t}\cdot\nabla\bar{u}_m \right)dxdt=0 .
	\end{equation}
	Since
	\begin{equation}
		\left|\frac{\left|\nabla \bar{u}_m(x,t)\right|^2}{1+\left|\nabla \bar{u}_m(x,t)\right|^2} - \frac{\left|\nabla u(x,t)\right|^2}{1+\left|\nabla u(x,t)\right|^2}\right|\leq \left|\left|\nabla \bar{u}_m(x,t)\right|-\left|\nabla u(x,t)\right|\right|\leq \left|\nabla \bar{u}_m(x,t)-\nabla u(x,t)\right|, \notag
	\end{equation}
	for a.e. $(x,t)\in Q_T$. Combining this with \eqref{strong convergen barum}, \eqref{Lambda equal um sigma}, and \eqref{Lambda equal um px}, we deduce that
	\begin{align}
		\iint_{Q_T}\delta\vec{\sigma}\cdot\nabla u dxdt = & - \iint_{Q_T} \frac{\partial u}{\partial t} u dxdt - \iint_{Q_T} \frac{\left|\nabla u\right|^2}{1+\left|\nabla u\right|^2} dxdt - \iint_{Q_T} \varepsilon \nabla \frac{\partial u}{\partial t}\cdot\nabla u dxdt \notag\\
		= & \lim_{m\rightarrow \infty} \left(-\iint_{Q_T}\frac{\partial u_m}{\partial t}\bar{u}_m dxdt - \iint_{Q_T} \frac{\left|\nabla \bar{u}_m\right|^2}{1+\left|\nabla \bar{u}_m\right|^2} dxdt - \iint_{Q_T} \varepsilon \nabla \frac{\partial u_m}{\partial t}\cdot\nabla\bar{u}_m dxdt\right) \notag\\
		= & \limsup_{m\rightarrow \infty} \iint_{Q_T}\delta\left|\nabla \bar{u}_m\right|^{p(x)} dxdt . \notag
	\end{align}
	Given that the mapping $\xi\mapsto \frac{1}{p}\left|\xi\right|^{p}$ is convex for $p\geq 1$, it follows that for any $\psi\in L^{\infty}(0,T;W^{1,p(\cdot)}(\Omega))$, we have
	\begin{equation}
		\iint_{Q_T}\left( \left|\nabla \bar{u}_m\right|^{p(x)-2}\nabla \bar{u}_m - \left|\nabla \psi\right|^{p(x)-2}\nabla \psi \right)\cdot\left(\nabla \bar{u}_m - \nabla \psi\right)dxdt \geq 0 . \notag
	\end{equation}
	Passing to the limit as $m \rightarrow \infty$ in the above inequality, we derive
	\begin{equation}
		\iint_{Q_T}\left(\vec{\sigma} - \left|\nabla \psi\right|^{p(x)-2}\nabla \psi\right)\cdot\left(\nabla u - \nabla \psi\right) dxdt \geq 0 .\notag
	\end{equation}
	Thus, applying the standard monotone method, we obtain $\vec{\sigma}=\left|\nabla u\right|^{p(x)-2}\nabla u$ a.e. in $Q_T$, which verifies that \eqref{Lambda equal} holds.
\end{proof}

Based on the results established earlier, we now prove that the regularized problem \eqref{P4}--\eqref{P6} admits a unique weak solution.

\begin{proof}[Proof of Theorem \ref{thm 1}]
	Combining \eqref{eqn Lambda} and \eqref{Lambda equal} yields the following integral identity
	\begin{equation}
		\iint_{Q_T}\left(\frac{\partial u}{\partial t}\zeta + \vec{q}(x,\nabla u)\cdot \nabla\zeta + \varepsilon \nabla \frac{\partial u}{\partial t}\cdot\nabla\zeta \right)dxdt=0 , \label{eqn weak solution without varepsilon}
	\end{equation}
	for every $\zeta\in L^1(0,T;W^{1,p(\cdot)}(\Omega))\cap L^2(0,T;H^1(\Omega))$. The Neumann boundary condition of problem \eqref{P4}--\eqref{P6} is implicitly incorporated in the weak variational formulation $J^{m,j}$ and can be recovered via integration by parts. Next, we address the uniqueness of the weak solution to problem \eqref{P4}--\eqref{P6}. To this end, assume that $u_1$ and $u_2$ are two weak solutions with the same initial data $f$. For any $t\in [0,T]$, choosing the test function $\zeta = \left(u_1 - u_2\right)\chi_{[0,t]}(\tau)$ in \eqref{eqn weak solution without varepsilon}, we obtain
	\begin{align}
		&\iint_{Q_t}\left(\frac{\partial u_1}{\partial \tau}-\frac{\partial u_2}{\partial \tau}\right)\left(u_1 - u_2\right)dxd\tau \notag\\
		+ & \iint_{Q_t}\left(\frac{\nabla u_1}{1 + \left|\nabla u_1\right|^2} - \frac{\nabla u_2}{1 + \left|\nabla u_2\right|^2}\right)\cdot\left(\nabla u_1 - \nabla u_2\right) dxd\tau \notag \\
		+ & \delta\iint_{Q_t}\left(\left|\nabla u_1\right|^{p(x)-2}\nabla u_1 - \left|\nabla u_2\right|^{p(x)-2}\nabla u_2\right)\cdot\left(\nabla u_1 - \nabla u_2\right)dxd\tau \notag\\
		+ & \varepsilon\iint_{Q_t}\left(\nabla \frac{\partial u_1}{\partial \tau} - \nabla \frac{\partial u_2}{\partial \tau}\right)\cdot\left(\nabla u_1 - \nabla u_2\right) dxd\tau = 0. \notag
	\end{align}
	The global Lipschitz continuity of the mapping $\xi\mapsto \frac{\xi}{1+\left|\xi\right|^2}$ implies that
	\begin{align}
		& \frac{1}{2}\int_{\Omega} \left|u_1(t) - u_2(t)\right|^2 dx + \frac{\varepsilon}{2}\int_{\Omega}\left|\nabla u_1(t) - \nabla u_2(t)\right|^2 dx\notag\\
		\leq & \iint_{Q_t}\left(\frac{\nabla u_1}{1 + \left|\nabla u_1\right|^2} - \frac{\nabla u_2}{1 + \left|\nabla u_2\right|^2}\right)\cdot\left(\nabla u_2 - \nabla u_1\right) dxd\tau \notag\\
		\leq & K(N) \int_{0}^{t}\int_{\Omega} \left|\nabla u_1 - \nabla u_2\right|^2 dxd\tau . \notag
	\end{align}
	Using Gronwall's inequality, we derive that $u_1 = u_2$ a.e. in $Q_T$, which shows the uniqueness of the weak solution to problem \eqref{P4}--\eqref{P6}.

	To prove Theorem \ref{thm 1}, it remains to establish the energy estimate \eqref{solution estimate}. By Lemma \ref{discrete energy inequality lemma}, for any $\gamma\in (0,\varepsilon)$, there exists $m_0(\gamma)>0$ such that, for all $m\geq m_0$ and for a.e. $t\in [0,T)$, we derive from \eqref{discrete energy inequality} the following inequalities, where $(j-1)\frac{T}{m}\leq t<j\frac{T}{m}$ holds for some $j$:
	\begin{align}
		&\int_{\Omega}\varphi(x,\nabla \bar{u}_m(t)) dx + \int_{0}^{t} \int_{\Omega} \left|\frac{\partial u_m}{\partial \tau}\right|^2 dxd\tau + \int_{0}^{t}\int_{\Omega}\left(\varepsilon - \gamma\right)\left|\nabla \frac{\partial u_m}{\partial \tau}\right|^2 dxd\tau \notag\\
		\leq & \int_{\Omega}\varphi(x,\nabla u_m^j) dx + \sum_{i=1}^{j}\int_{(i-1)\frac{T}{m}}^{i\frac{T}{m}} \int_{\Omega} \left|\frac{\partial u_m}{\partial t}\right|^2 dxdt + \sum_{i=1}^{j}\int_{(i-1)\frac{T}{m}}^{i\frac{T}{m}}\int_{\Omega}\left(\varepsilon - \gamma\right)\left|\nabla \frac{\partial u_m}{\partial t}\right|^2 dxdt \notag\\
		= & \int_{\Omega}\varphi(x,\nabla u_m^j) dx + \sum_{i=1}^{j} \int_{\Omega} \frac{m}{T} \left|u_m^i-u_m^{i-1}\right|^2 dx + \sum_{i=1}^{j} \int_{\Omega} \left(\varepsilon - \gamma\right) \frac{m}{T} \left|\nabla u_m^i-\nabla u_m^{i-1}\right|^2 dx \notag\\
		\leq & \int_{\Omega}\varphi(x,\nabla f) dx . \notag
	\end{align}
	Since $\frac{\partial u_m}{\partial \tau}\rightharpoonup \frac{\partial u}{\partial \tau}$ weakly in $L^2(0,t;H^1(\Omega))$, it follows from the weak lower semicontinuity of the $L^2(Q_t)$-norm that
	\begin{equation}
		\liminf_{m\rightarrow \infty}\int_{\Omega}\varphi(x,\nabla \bar{u}_m(t)) dx + \int_{0}^{t} \int_{\Omega} \left|\frac{\partial u}{\partial \tau}\right|^2 dxd\tau + \int_{0}^{t}\int_{\Omega}\left(\varepsilon - \gamma\right)\left|\nabla \frac{\partial u}{\partial \tau}\right|^2 dxd\tau \leq \int_{\Omega}\varphi(x,\nabla f) dx . \notag
	\end{equation}
	From Lemma \ref{strong convergence gradient um}, we know that $\nabla \bar{u}_m(t)\rightarrow\nabla u(t)$ a.e. in $\Omega$. Then, Fatou's Lemma implies that
	\begin{equation}
		\int_{\Omega}\varphi(x,\nabla u(t)) dx + \int_{0}^{t} \int_{\Omega} \left|\frac{\partial u}{\partial \tau}\right|^2 dxd\tau + \int_{0}^{t}\int_{\Omega}\left(\varepsilon - \gamma\right)\left|\nabla \frac{\partial u}{\partial \tau}\right|^2 dxd\tau \leq \int_{\Omega}\varphi(x,\nabla f) dx . \notag
	\end{equation}
	Given that $\gamma > 0$ is arbitrary, we obtain
	\begin{equation}
		\int_{\Omega}\varphi(x,\nabla u(t)) dx + \int_{0}^{t} \int_{\Omega} \left|\frac{\partial u}{\partial \tau}\right|^2 dxd\tau + \int_{0}^{t}\int_{\Omega}\varepsilon\left|\nabla \frac{\partial u}{\partial \tau}\right|^2 dxd\tau \leq \int_{\Omega}\varphi(x,\nabla f) dx , \notag
	\end{equation}
	for a.e. $t\in [0,T]$. This completes the proof of Theorem \ref{thm 1}.
\end{proof}

\section{Proof of Theorem \ref{main thm}}\label{Proof of Theorem 1}

In this section, we prove the existence of Young measure solutions to problem \eqref{P1}--\eqref{P3} using Theorem \ref{thm 1} and the corresponding limit process $\varepsilon \to 0^{+}$.

For any $\varepsilon \in (0,1)$, it follows from Theorem \ref{thm 1} that the regularized problem \eqref{P4}--\eqref{P6} admits a unique weak solution $u^{\varepsilon}\in L^{\infty}(0,T;W_f^{1,p(\cdot)}(\Omega))\cap H^{1}(0,T;H^1(\Omega))$, which satisfies
\begin{equation}
	\iint_{Q_T}\left(\frac{\partial u^{\varepsilon}}{\partial t}\zeta + \vec{q}(x,\nabla u^{\varepsilon})\cdot \nabla\zeta + \varepsilon \nabla \frac{\partial u^{\varepsilon}}{\partial t}\cdot\nabla\zeta \right)dxdt=0 , \label{eqn:3}
\end{equation}
for every $\zeta\in L^1(0,T;W^{1,p(\cdot)}(\Omega))\cap L^2(0,T;H^1(\Omega))$, together with the energy estimate
\begin{equation}
	\int_{\Omega}\varphi(x,\nabla u^{\varepsilon}(t))dx+\left\| \frac{\partial u^{\varepsilon}}{\partial \tau} \right\|^2_{L^2(Q_t)} + \left\| \nabla \sqrt{\varepsilon} \frac{\partial u^{\varepsilon}}{\partial \tau} \right\|^2_{L^2(Q_t)} \leq \int_{\Omega}\varphi(x,\nabla f) dx , \notag
\end{equation}
for a.e. $t\in [0,T]$. Then, from \eqref{structure condition 1}, we derive
\begin{equation}
	\varrho_{L^{p(\cdot)}(\Omega)}\left(\nabla u^{\varepsilon}(t)\right) \leq \lambda_1^{-1} \left(\int_{\Omega}\varphi(x,\nabla u^{\varepsilon}(t)) dx + \left|\Omega\right|\right) \leq \lambda_1^{-1} \left(\int_{\Omega}\varphi(x,\nabla f) dx + \left|\Omega\right|\right) , \notag
\end{equation}
and furthermore,
\begin{align}
	\varrho_{L^{p^{\prime}(\cdot)}(\Omega)}\left(\vec{q}(x,\nabla u^{\varepsilon}(t))\right) = & \int_{\Omega} \left|\vec{q}(x,\nabla u^{\varepsilon}(t))\right|^{\frac{p(x)}{p(x)-1}} dx \notag\\
	\leq & \int_{\Omega} 2^{\frac{1}{p(x)-1}}\left( \lambda_{2}^{\frac{p(x)}{p(x)-1}}\left|\nabla u^{\varepsilon}(t)\right|^{p(x)} + 1\right) dx \notag\\
	\leq & 2^{\frac{1}{p^{-}-1}} \left(\lambda_2+1\right)^{\frac{p^{-}}{p^{-}-1}} \varrho_{L^{p(\cdot)}(\Omega)}\left(\nabla u^{\varepsilon} (t)\right) + 2^{\frac{1}{p^{-}-1}}\left|\Omega\right| \notag\\
	\leq & 2^{\frac{1}{p^{-}-1}} \left(\lambda_2+1\right)^{\frac{p^{-}}{p^{-}-1}} \lambda_1^{-1} \left(\int_{\Omega}\varphi(x,\nabla f) dx + \left|\Omega\right|\right) + 2^{\frac{1}{p^{-}-1}}\left|\Omega\right| .  \notag
\end{align}
By Lemmas \ref{norm vs semi norm} and \ref{Poincare inequality general trace}, along with the mean invariance of $u^{\varepsilon}(t)$, we deduce that
\begin{equation}
	\left\|u^{\varepsilon}_t\right\|_{L^2(Q_T)} + 	\left\|\nabla \sqrt{\varepsilon} u^{\varepsilon}_t\right\|_{L^2(Q_T)} + \sup_{0\leq t\leq T} \left\| u^{\varepsilon}(t) \right\|_{W^{1,p(\cdot)}(\Omega)}
	+ \sup_{0\leq t\leq T} \left\| \vec{q}\left(x, \nabla u^{\varepsilon}(t)\right) \right\|_{L^{p^{\prime}(\cdot)}(\Omega)} \leq M , \label{estimate grau eps}
\end{equation}
where $M>0$ is independent of $\varepsilon$. Therefore, there exist a subsequence of $\{u^{\varepsilon}\}_{0<\varepsilon<1}$, still denoted by itself, and a function $u \in L^{\infty}(0,T;W_f^{1,p(\cdot)}(\Omega))\cap H^1(0,T;L^2(\Omega))$ such that the following weak convergence results hold as $\varepsilon \to 0^{+}$:
\begin{align}
	& u^{\varepsilon} \rightharpoonup u \text{ weakly }{\ast}\text{ in } L^{\infty}(0,T;W^{1,p(\cdot)}(\Omega)), \notag\\
	& u^{\varepsilon} \rightarrow u \text{ strongly in } C([0,T];L^{2}(\Omega)), \notag\\
	& u^{\varepsilon}_t \rightharpoonup u_t \text{ weakly in } L^2(Q_T) . \notag
\end{align}

For a.e. $t \in [0,T]$, let $\nu_{\cdot,t} = (\nu_{x,t})_{x \in \Omega}$ denote the Young measure generated by the sequence $\{\nabla u^{\varepsilon}(t)\}_{0<\varepsilon<1}$. Since the generating sequence satisfies the boundedness condition \eqref{boundedness condition}, it follows from Lemmas \ref{Fundamental theorem on Young measures} and \ref{equivalent lemma} that $\left\|\nu_{x,t}\right\|_{\mathcal{M}(\mathbb{R}^N)} = 1$ a.e. in $\Omega$. We now verify that the pair $\left(u,\nu\right)$ is a Young measure solution to problem \eqref{P1}--\eqref{P3}. Combining Remark \ref{rmk 1} with \eqref{estimate grau eps}, we infer that
\begin{equation}
	\vec{q}(\cdot,\nabla u^{\varepsilon}(t)) \rightharpoonup \langle\nu_{x,t},\vec{q}(x,\cdot) \rangle=\int_{\mathbb{R}^N}\vec{q}(\cdot,\xi)d\nu_{\cdot,t}(\xi)\text{ weakly in } L^{p^{\prime}(\cdot)}(\Omega), \notag
\end{equation}
which yields
\begin{equation}
	\int_{\Omega}\vec{q}(x,\nabla u^{\varepsilon}(t)) \cdot \vec{\psi}(t) dx \rightarrow \int_{\Omega}\langle\nu_{x,t},\vec{q}(x,\cdot)\rangle \cdot \vec{\psi}(t) dx\text{ a.e. in }[0,T], \notag
\end{equation}
for every $\vec{\psi}\in L^1(0,T;L^{p(\cdot)}(\Omega))$. By the dominated convergence theorem, we obtain
\begin{equation}
	\int_{0}^{T}\int_{\Omega}\vec{q}(x,\nabla u^{\varepsilon}) \cdot \vec{\psi} dxdt \rightarrow \int_{0}^{T}\int_{\Omega}\langle\nu_{x,t},\vec{q}(x,\cdot)\rangle \cdot \vec{\psi} dxdt, \notag
\end{equation}
and thus
\begin{equation}
	\vec{q}(\cdot,\nabla u^{\varepsilon}) \rightharpoonup \langle\nu,\vec{q}(x,\cdot)\rangle\text{ weakly }{\ast}\text{ in }L^{\infty}(0,T;L^{p^{\prime}(\cdot)}(\Omega)) . \notag
\end{equation}
In a similar manner, we have
\begin{equation}
	\nabla u^{\varepsilon} \rightharpoonup \langle\nu,\mathrm{id}\rangle \text{ weakly }{\ast}\text{ in }L^{\infty}(0,T;L^{p(\cdot)}(\Omega)). \notag
\end{equation}
Then, we obtain $\nabla u (x,t) = \langle\nu_{x,t},\mathrm{id}\rangle$ a.e. in $Q_T$, proving \eqref{gradient}. Furthermore, using H\"{o}lder's inequality, we derive
\begin{align}
	\left|\iint_{Q_T}\varepsilon \nabla \frac{\partial u^{\varepsilon}}{\partial t}\cdot\nabla\zeta dxdt\right| \leq \sqrt{\varepsilon}\left\|\nabla\sqrt{\varepsilon}\frac{\partial u^{\varepsilon}}{\partial t} \right\|_{L^2(Q_T)} \left\|\nabla\zeta \right\|_{L^2(Q_T)} \leq \sqrt{\varepsilon M} \left\|\nabla\zeta \right\|_{L^2(Q_T)}. \notag
\end{align}
Taking the limit $\varepsilon\rightarrow 0^{+}$ in \eqref{eqn:3}, we arrive at
\begin{equation}
	\iint_{Q_T}\left(\frac{\partial u}{\partial t}\zeta + \langle\nu_{x,t},\vec{q}(x,\cdot) \rangle\cdot \nabla\zeta\right)dxdt=0 , \notag
\end{equation}
for every $\zeta\in L^1(0,T;W^{1,p(\cdot)}(\Omega))\cap L^2(0,T;H^1(\Omega))$. Moreover, $u(0)=f$ a.e. in $\Omega$ follows from the same argument as in the proof of Lemma \ref{strong convergence gradient um}, and the details are omitted here. Therefore, $\left(u,\nu\right)$ is a Young measure solution to problem \eqref{P1}--\eqref{P3}, and the proof is complete.

% ------------------------------------------------------------------------

\subsection*{Acknowledgment}
The authors wish to express their sincere gratitude to Lecturer Jingfeng Shao (Guangxi University) for many helpful discussions during the preparation of this paper. This work was supported by the National Natural Science Foundation of China (Grant Nos. U21B2075, 12171123, 12301536) and the Fundamental Research Funds for the Central Universities (Grant Nos. 2022FRFK060020, HIT.NSRIF202302).

\section*{\small
 Conflict of interest} %%%%%%%%%%%%%%%%%%%%

 {\small
 The authors declare that they have no conflict of interest.}

% ------------------------------------------------------------------------
\end{document}